\documentclass[11pt,a4paper]{article}

\usepackage[utf8]{inputenc}
\usepackage[T1]{fontenc}
\usepackage{tikz}
\usepackage{pgfplots}
\pgfplotsset{compat=1.16}

\usepackage[a4paper,margin=1in]{geometry}
\usepackage{setspace}
\usepackage{amsmath,amssymb,amsfonts,amsthm}
\usepackage{mathtools}
\usepackage{bm}
\usepackage{comment}

\usepackage{graphicx}
\usepackage{subcaption}
\usepackage{booktabs}
\usepackage{array}
\usepackage{float}
\usepackage{placeins}
\usepackage{flafter}
\usepackage{enumitem}
\usepackage{multirow}
\usepackage{longtable}

\usepackage{algorithm}
\usepackage{algpseudocode}

\usepackage[numbers,sort&compress]{natbib}
\usepackage{xcolor}
\usepackage[
    colorlinks=true,
    linkcolor=blue,
    citecolor=blue,
    urlcolor=blue
]{hyperref}
\usepackage{aliascnt}
\usepackage[nameinlink,capitalise,noabbrev]{cleveref}

\theoremstyle{plain}
\newtheorem{theorem}{Theorem}[section]
\newaliascnt{proposition}{theorem}
\newtheorem{proposition}[proposition]{Proposition}
\aliascntresetthe{proposition}
\crefname{proposition}{Proposition}{Propositions}
\Crefname{proposition}{Proposition}{Propositions}
\newaliascnt{lemma}{theorem}
\newtheorem{lemma}[lemma]{Lemma}
\aliascntresetthe{lemma}
\crefname{lemma}{Lemma}{Lemmas}
\Crefname{lemma}{Lemma}{Lemmas}
\newaliascnt{corollary}{theorem}
\newtheorem{corollary}[corollary]{Corollary}
\aliascntresetthe{corollary}
\crefname{corollary}{Corollary}{Corollaries}
\Crefname{corollary}{Corollary}{Corollaries}

\theoremstyle{definition}
\newaliascnt{definition}{theorem}
\newtheorem{definition}[definition]{Definition}
\aliascntresetthe{definition}
\crefname{definition}{Definition}{Definitions}
\Crefname{definition}{Definition}{Definitions}
\newaliascnt{assumption}{theorem}
\newtheorem{assumption}[assumption]{Assumption}
\aliascntresetthe{assumption}
\crefname{assumption}{Assumption}{Assumptions}
\Crefname{assumption}{Assumption}{Assumptions}
\newaliascnt{example}{theorem}

\aliascntresetthe{example}
\crefname{example}{Example}{Examples}
\Crefname{example}{Example}{Examples}

\theoremstyle{remark}
\newaliascnt{remark}{theorem}

\aliascntresetthe{remark}
\crefname{remark}{Remark}{Remarks}
\Crefname{remark}{Remark}{Remarks}

\newcommand{\R}{\mathbb{R}}

\newcommand{\norm}[1]{\left\lVert #1 \right\rVert}

\newcommand{\argmin}{\operatorname*{arg\,min}}
\newcommand{\diag}{\operatorname{diag}}

\title{
    A Structured-Surrogate Double-Model Approach\\
    to Regularization Parameter Selection
}

\author{
Malena Español \thanks{School of Mathematical and Statistical Sciences, Arizona State University, 1151 S Forest Ave, 85287 Tempe, AZ, United States, email: \texttt{mespanol@asu.edu}} 
\and
Markus Juvonen \thanks{Department of Mathematics and Statistics, University of Helsinki, Pietari Kalmin katu 5, 00560 Helsinki, Finland, email: \texttt{markus.juvonen@helsinki.fi}} 
\and
Samuli Siltanen \thanks{Department of Mathematics and Statistics, University of Helsinki, Pietari Kalmin katu 5, 00560 Helsinki, Finland, email: \texttt{samuli.siltanen@helsinki.fi}}
\and
Ilaria Trombini \thanks{Department of Mathematics and Computer Science, University of Ferrara, Via Machiavelli 30, 44121 Ferrara, Italy,  email: \texttt{ilaria.trombini@unife.it} } 
\and
Chuyang Wu \thanks{Department of Mathematics and Statistics, University of Helsinki, Pietari Kalmin katu 5, 00560 Helsinki, Finland, email: \texttt{chuyang.wu@helsinki.fi}} 
} 

\date{}

\begin{document}

\maketitle

\begin{abstract}
\label{sec:abstract}
Double-model methods select a regularization parameter by matching the discrepancy between two reconstructions of the same data to a prescribed target. The Structured-Surrogate Double-Model (SSDM) framework is developed around a computationally inexpensive structured surrogate of the forward operator and a spectrally explicit perturbation, with the selected parameter subsequently applied to the original inverse problem. Given a prescribed target, parameter selection requires neither prior knowledge of the noise level nor training data. For a class of Tikhonov regularization problems, the squared discrepancy is shown to attain a strict maximum at zero regularization, decay as $\lambda^{-2}$, and decrease strictly beyond an explicit phase-aware threshold. A complementary, data-dependent certificate is established to guarantee strict monotonicity on finite intervals, enabling unique and reliable target matching even when the global threshold is overly conservative. The discrepancy is further connected to directional roughness in image deblurring, providing a geometric interpretation of the selection principle. Experiments on synthetic images and on synthetically blurred natural photographs illustrate the theoretical findings and demonstrate the practical applicability of certified parameter selection.
\end{abstract}

\section{Introduction}
\label{sec:introduction}

We consider discrete linear inverse problems
\begin{equation}
b=Ax+\eta,
\label{eq:intro-forward-model}
\end{equation}
where $A\in\R^{m\times n}$ is a known forward operator, $b\in\R^m$ is
observed data, and $\eta$ is unknown measurement noise. Such problems are
typically ill-posed: $A$ is ill-conditioned or singular, and the naive
least-squares solution amplifies the noise far beyond any usable signal.
Tikhonov regularization~\cite{tikhonov1963solution, hoerl1970ridge} replaces it by
\begin{equation}
x_\lambda\in\argmin_x
\left\{\|Ax-b\|_2^2+\lambda\|Lx\|_2^2\right\},
\label{eq:intro-tikhonov}
\end{equation}
where $L$ encodes prior information and the penalty weight $\lambda>0$
balances data fit against the penalty. The penalty is often written as $\alpha^2\|Lx\|_2^2$. In this work, we use $\lambda$, where $\lambda=\alpha^2$.

Choosing $\lambda$ is a challenge. Established parameter-choice approaches require information that may be unavailable or difficult to estimate. For instance, the discrepancy principle requires the noise level~\cite{morozov1966solution}, and generalized cross-validation requires traces or spectral information of $A$~\cite{golub1979generalized}. The L-curve avoids both requirements but relies on a well-defined corner, which may be ambiguous or absent~\cite{hansen1992analysis}. Learning-based approaches need representative training data and, in supervised formulations, ground-truth solutions~\cite{chung2017learning,kunisch2013,arridge2019}. For large problems with imperfectly known physics, these requirements are often unavailable, expensive, or unreliable.

Two ideas offer a way around these requirements and difficulties. The first is to let a cheaper model serve as a proxy for $A$ when choosing $\lambda$. In image deblurring, $A$ blurs a scene that continues beyond the field of view. A \emph{structured surrogate} $\widetilde A$ that imposes periodic boundary conditions is inaccurate near the boundary, but it is diagonalized by the Discrete Fourier Transform (DFT), so every Tikhonov reconstruction with $\widetilde A$ costs two transforms, which can be computed efficiently via the Fast Fourier Transform (FFT). This is the idea behind the operator-approximation framework introduced in~\cite{chung2015framework}, which chooses the parameter $\lambda$ using an approximate operator and then transfers it to the original problem. 

The second idea is to reconstruct the same data twice using two deliberately different models and measure how much the reconstructions disagree~\cite{niinimaki2016,juvonen2025,wu2025}. This idea is explored in the dual-grid method presented in~\cite{juvonen2025}, where the second model is the first model applied on a grid shifted by half a pixel. When $\lambda$ is small, both reconstructions are dominated by amplified noise and by detail that depends on the precise model, and they disagree. As $\lambda$ grows, both reconstructions are smoothed, and they come to agree. Selecting the $\lambda$ at which the disagreement reaches a prescribed level $\tau$ gives a parameter choice that needs neither a noise level nor training data. In the context of tomography, the results in~\cite{wu2025} show that a proportional feedback controller can drive the disagreement to $\tau$. Related feedback principles appear in the controlled-sparsity approaches of~\cite{Purisha_2018,purisha2018automatic,Meaney2025image}.

We propose a double-model rule that combines both ideas. For each $\lambda$, it computes two Tikhonov reconstructions from two different surrogate models and measures how much they disagree. The rule selects the $\lambda$ at which this disagreement equals an acceptable level, the target~$\tau$.  The target is taken as prescribed: it could be supplied by the application, calibrated using representative data, learned, or set based on user experience. The rule is well-defined only if matching the target determines a parameter, which requires that the discrepancy between the two reconstructions decrease as $\lambda$ grows. Earlier double-model work assumes this or observes it numerically when using an image-similarity (SSIM) measure to define the disagreement, leaving uniqueness and formal conditions as open research questions~\cite{juvonen2025,wu2025}. For the classical discrepancy principle, the corresponding property is automatic, because the residual of the Tikhonov solution is monotone in $\lambda$~\cite{engl1996}. For the double-model discrepancy, it is not. Even in one dimension, a perturbation that merely rescales the model makes the discrepancy dip to zero, rise, and then decay, crossing the target several times (\cref{app:controller-multiple-crossings}). We therefore ask when a target has a unique solution, and how to compute it.

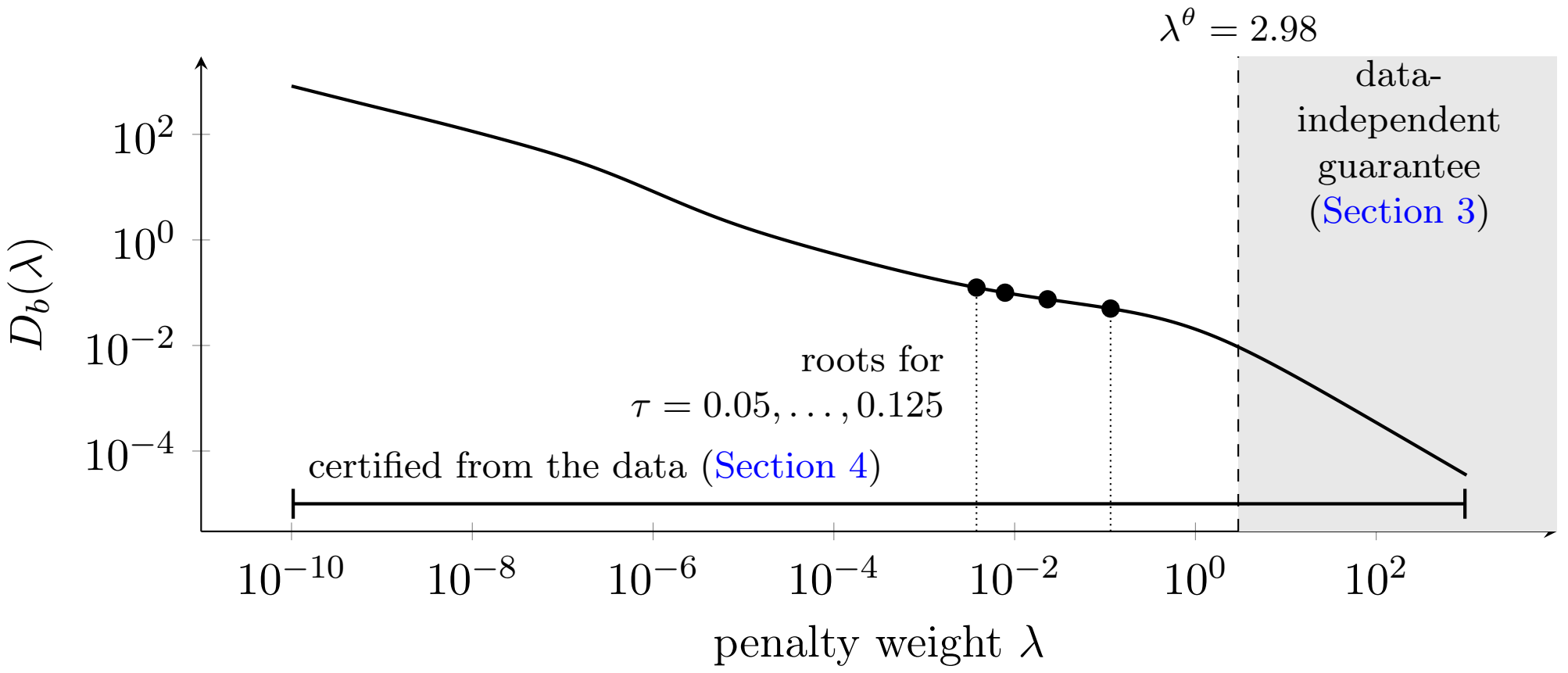
\begin{figure}[tb]
\centering
\begin{tikzpicture}
\begin{axis}[
  width=0.82\textwidth, height=5.6cm,
  xmode=log, ymode=log,
  xmin=1e-11, xmax=1e4, ymin=3e-6, ymax=3e3,
  xlabel={penalty weight $\lambda$},
  ylabel={$D_b(\lambda)$},
  xtick={1e-10,1e-8,1e-6,1e-4,1e-2,1,1e2},
  ytick={1e-4,1e-2,1,1e2},
  axis lines=left, clip=false,
  every axis plot/.append style={thick}]
  % data-independent tail of Theorem 3.8: (lambda^theta, infinity)
  \fill[gray!18] (axis cs:2.977223661,3e-6) rectangle (axis cs:1e4,3e3);
  \node[align=center,font=\footnotesize] at (axis cs:180,60)
    {data-\\independent\\guarantee\\(\cref{sec:spectral-theory})};
  \draw[dashed] (axis cs:2.977223661,3e-6) -- (axis cs:2.977223661,3e3)
    node[above,font=\footnotesize] {$\lambda^\theta=2.98$};
  % D_b for Experiment 2, case Z-B, at every 10th certificate cell edge
  \addplot[black] table[x=lambda, y=Db, col sep=comma] {fig1_curve.csv};
  % roots for tau = 0.05, 0.075, 0.10, 0.125
  \addplot[only marks, mark=*, mark size=1.8pt, black] coordinates {
    (0.1151774776,0.050) (0.02312570272,0.075)
    (0.007833447765,0.100) (0.003783693796,0.125)};
  \draw[densely dotted] (axis cs:0.003783693796,3e-6) -- (axis cs:0.003783693796,0.125);
  \draw[densely dotted] (axis cs:0.1151774776,3e-6) -- (axis cs:0.1151774776,0.050);
  \node[font=\footnotesize,anchor=north east,align=right] at (axis cs:2.5e-3,0.02)
    {roots for\\$\tau=0.05,\ldots,0.125$};
  % certified interval [1e-10, 1e3]
  \draw[|-|, thick] (axis cs:1e-10,1e-5) -- (axis cs:1e3,1e-5)
    node[pos=0,anchor=south west,font=\footnotesize]
    {certified from the data (\cref{sec:controller})};
\end{axis}
\end{tikzpicture}
\caption{The two levels of the analysis, on the zero-boundary data of
Experiment~2 (bilinear perturbation, periodic surrogate). The spectral
theory of \cref{sec:spectral-theory} guarantees that the normalized
disagreement $D_b$ decreases beyond the threshold $\lambda^\theta=2.977$. The roots of $D_b=\tau$ for the four prescribed targets (dots) lie $1.4$ to $2.9$ decades below $\lambda^\theta$. There, decrease is a property of the data: the certificate of
\cref{sec:controller} verifies it on the whole interval $[10^{-10},10^3]$, so each target is met at a unique root.}
\label{fig:two-levels}
\end{figure}

To address these questions, we introduce the \emph{Structured-Surrogate Double-Model} (SSDM) framework that can be summarized as follows:
\begin{equation}
A
\;\longrightarrow\;
\underbrace{\widetilde A}_{\text{structured surrogate}}
\;\longrightarrow\;
\underbrace{(\widetilde A,\widetilde AS)}_{\text{double model}}
\;\overset{\text{target }\tau}{\longrightarrow}\;
\lambda_{\mathrm{SSDM}}
\;\longrightarrow\;
\underbrace{x_{\lambda_{\mathrm{SSDM}}}(A)}_{\text{original problem}}.
\label{eq:ssdm-pipeline}
\end{equation}
Its defining choice is that the second model $\widetilde AS$ is obtained from the surrogate by a perturbation $S$ that is diagonal in the surrogate's own basis. Both reconstructions are then filtered mode-wise, and the squared disagreement becomes a weighted sum of explicit one-mode functions whose weights are the data energies. This makes the disagreement analyzable, and the analysis takes place at two levels (\cref{fig:two-levels}).

At the \emph{data-independent} level, for a surrogate, its perturbation, and a regularization operator that are simultaneously unitarily diagonalizable, and a perturbation that amplifies no mode, the disagreement has a clean global shape. It is largest at zero regularization, decays like $\lambda^{-1}$, and decreases strictly beyond a threshold computed from the spectra alone. Tracking the relative phase between the two models sharpens this threshold optimally: under a mild sign-change condition, no argument that ignores the data can establish a decrease below it. This motivates the \emph{data-dependent} level. In our half-pixel deblurring example, the threshold lies between about $1.6$ and $3$, whereas the parameters selected experimentally lie between $3\times10^{-4}$ and $0.12$. In the practical range, monotonicity is a property of the data. We give a finite certificate: on the cells of a partition, it bounds the derivative of the disagreement from above, using one cubic polynomial per mode. Where the certificate holds, every target in range is met at exactly one parameter, and a bracketing solver finds it in about a dozen evaluations.

For the half-pixel pair, the disagreement also has a concrete meaning. It is a data-normalized directional roughness of the surrogate reconstruction: exactly so for a pure-phase half-pixel shift, and to leading order at low frequency for bilinear interpolation. A target, therefore, prescribes how much variation along the shift direction the reconstruction may carry, relative to the data, a quantity that can be reasoned about directly.

Three deblurring experiments put the theory to work. The first verifies it
in a matched model. The second transfers the selected parameter to a blur
with zero boundary values that the periodic surrogate does not reproduce.
There, rules that match the surrogate residual to the noise level fail even
when given the true noise norm, and recover only when the level is inflated
to absorb the model error, while the SSDM disagreement remains usable
without forming that residual. The third blurs three bird photographs on the computer, so that the ground truth is known, and reconstructs them with a blur model that is slightly wrong and a surrogate whose boundary behavior is also wrong. In every configuration where it was computed, floating-point evaluation of the certificate verifies its conditions over the whole searched range. 

In summary, our contributions are:
\begin{itemize}[leftmargin=*]
\item the global shape of the SSDM disagreement, with a phase-aware
monotonicity threshold that is optimal among data-independent thresholds
(\cref{thm:main-spectral-general,thm:phase-aware-threshold});
\item a finite, data-dependent certificate of strict decrease, which makes
target matching well posed in the range used in practice, and a simple
certified selection procedure (\cref{prop:interval-certificate},
\cref{alg:ssdm});
\item for half-pixel deblurring, explicit thresholds and the interpretation
of the disagreement as directional roughness
(\cref{cor:half-pixel-thresholds,prop:roughness});
\item experiments on synthetic data and on natural photographs that verify the certificate,
examine parameter transfer under boundary mismatch, and show how the target
behaves (\cref{sec:numerical-experiments}).
\end{itemize}

Our analysis is specific to Tikhonov regularization with operators that
share an eigenbasis. However, the SSDM architecture is not limited to this specification: double-model selection
has been used with total-variation regularization~\cite{juvonen2025,wu2025}, and the framework of~\cite{chung2015framework} accommodates other spectral filters. We therefore regard
the present results as a first analytically tractable class within that
broader framework. 

The rest of the paper is organized as follows. \Cref{sec:framework} defines the SSDM framework. \Cref{sec:spectral-theory} develops the data-independent theory and \cref{sec:controller} the
certificate. \Cref{sec:problem-setup} specializes both to half-pixel
deblurring, \cref{sec:numerical-experiments} reports the experiments, and
\cref{sec:conclusion} concludes. Proofs that would interrupt the argument are in the appendices.

\section{The SSDM framework}
\label{sec:framework}

\paragraph{Tikhonov reconstructions.}
For a forward matrix $B\in\R^{m\times n}$, data $b\in\R^m$, a penalty
matrix $L\in\R^{r\times n}$, and $\lambda>0$, we write $x_\lambda(B)$ for
the minimum-norm minimizer of $\|Bx-b\|_2^2+\lambda\|Lx\|_2^2$,
\begin{equation}
x_\lambda(B)=(B^*B+\lambda L^*L)^\dagger B^*b,
\label{eq:minimum-norm-convention}
\end{equation}
where $\dagger$ is the Moore--Penrose pseudoinverse. If
$\ker B\cap\ker L=\{0\}$, the minimizer is unique and $\dagger$ is the
ordinary inverse. In our spectral setting, the assumption that the
penalty is positive on every mode excited by the data confines possible
nonuniqueness to inactive modes
(\cref{subsec:general-phase-formula}).

The SSDM pipeline \eqref{eq:ssdm-pipeline} has four steps. We describe each
with the running example that motivates the paper, image deblurring with a
half-pixel shift.

\paragraph{Step 1: a structured surrogate.}
For the spectral realization studied here, we restrict to
$A\in\R^{n\times n}$, nonzero data $b\in\R^n$, and
$L\in\R^{n\times n}$. Following~\cite{chung2015framework}, we replace $A$
by a structured surrogate $\widetilde A\in\R^{n\times n}$ that
approximates it in the same data and solution coordinates and is
diagonalized by a known unitary transform $Q$. The surrogate, the
perturbation of Step~2, and the penalty share this basis:
\begin{equation}
\widetilde A=Q\diag(\nu)Q^*,\qquad
S=Q\diag(\gamma)Q^*,\qquad
L=Q\diag(\nu^L)Q^*,
\label{eq:simultaneous-diagonalization}
\end{equation}
with complex eigenvalues $\nu_k$, $\gamma_k$, and $\nu^L_k$. The operators
and data are real. The complex basis is only a representation. The original
operator $A$ need not be diagonalized by $Q$. For a spatially invariant
blur, periodic boundary conditions give a BCCB surrogate diagonalized by
the two-dimensional discrete Fourier transform, and reflexive boundary
conditions with a symmetric point-spread function give one diagonalized
by the discrete cosine transform. A blur with zero boundary values is not
Fourier-diagonalizable, but it can be approximated systematically, for
example by the nearest BCCB matrix in the Frobenius
norm~\cite{chan1988optimal,chan1994circulant,kamm1998kronecker}. In each
case the surrogate spectrum costs a few fast transforms, even when the
spectrum of $A$ is inaccessible.

\paragraph{Step 2: a perturbed second model and its disagreement.}
The second model is obtained by composing the surrogate with a
perturbation $S$ diagonalized by the same $Q$:
$\widetilde A_S:=\widetilde AS$. In the running example $S$ is bilinear
half-pixel interpolation: it replaces each pixel by the average of a
$2\times2$ block, which is a convolution and hence Fourier-diagonal
(\cref{prop:half-pixel-multiplier}). With the same data, penalty, and
parameter for both models, let
$\widetilde x_\lambda:=x_\lambda(\widetilde A)$ and
$\widetilde x_\lambda^S:=x_\lambda(\widetilde A_S)$. Define their
discrepancy and its data-normalized form by
\begin{equation}
D(\lambda):=\|\widetilde x_\lambda-\widetilde x_\lambda^S\|_2,\qquad
D_b(\lambda):=\frac{D(\lambda)}{\|b\|_2},\qquad
\Phi:=D^2,\quad\Phi_b:=D_b^2.
\label{eq:framework-discrepancy}
\end{equation}
\label{sec:normalized-variant}%
The assumption $b\ne0$ makes $D_b$ well defined. The reconstructions are
linear in $b$, so $D_b$ is unchanged when $b$ is multiplied by a nonzero
scalar. Targets are set on this normalized scale. The perturbation should
make the two reconstructions differ where they remain sensitive to
model-specific detail, while regularization makes them agree.

% \begin{remark}[What the normalization does and does not remove]
% \label{rem:normalization-scope}
% Division by $\|b\|_2$ makes $D_b$ independent of the scale of the data, but
% not of everything that sets the scale of a target. (i)~Adding a constant to
% $b$ leaves the numerator unchanged whenever the zero-frequency multiplier of
% $S$ is $1$, because the constant mode then contributes nothing to $D$
% (\cref{prop:half-pixel-multiplier}), but it changes the normalizing data norm,
% $\|b+\kappa\mathbf1\|_2^2=\|b\|_2^2+2\kappa\mathbf1^\top b+n\kappa^2$, and
% therefore in general changes $D_b$.
% (ii)~$D_b$ carries the units of $x/b$, that is, of the inverse of the
% operator. Replacing $\widetilde A$ and $\widetilde A_S$ by $c\widetilde A$ and
% $c\widetilde A_S$, with $c>0$ and the penalty $L$ held fixed, changes
% $D_b(\lambda)$ to $c^{-1}D_b(\lambda/c^{2})$ (for real $c\ne0$ the prefactor
% is $|c|^{-1}$), so a target
% is comparable across problems only for operators of comparable gain. The original PSFs
% have unit mass. The same-PSF surrogates, therefore, have unit zero-frequency
% gain, while the optimal circulant surrogate, the Frobenius-nearest BCCB approximation of A~\cite{chan1988optimal}, has gain $0.98472$. A target is therefore a statement about data of a given offset and an operator of a given gain, as well as about the model pair. We have not tested whether these two effects
% contribute to the variation of the oracle-associated target in
% \cref{sec:target-choice}.
% \end{remark}

\paragraph{Step 3: match a prescribed target.}
For a target $\tau>0$, select $\lambda_{\mathrm{SSDM}}$ with
\begin{equation}
D_b(\lambda_{\mathrm{SSDM}})=\tau
\label{eq:framework-controller-target}
\end{equation}
to a prescribed tolerance. The target is an external input (\cref{sec:introduction}). \Cref{sec:spectral-theory,sec:controller} establish when \eqref{eq:framework-controller-target} has a unique solution on a known
interval and how to find it.

\paragraph{Step 4: transfer the parameter.}
Reconstruct with the original operator and the same penalty,
\begin{equation}
x_{\lambda_{\mathrm{SSDM}}}(A)
=(A^*A+\lambda_{\mathrm{SSDM}}L^*L)^\dagger A^*b .
\label{eq:framework-transfer}
\end{equation}

Steps~1 and~4 come from operator approximation, and Step~2 from the
double-model idea. What is new is that the two models are tied together by
a perturbation that is diagonal in the surrogate's own basis. The next
section shows what this buys us: the discrepancy decomposes into explicit
modal terms whose behavior can be analyzed exactly. The analysis does not
guarantee the quality of the transferred reconstruction, which depends
on the target and on how well $\widetilde A$ approximates $A$; these
effects are examined in the experiments.

\paragraph{Rectangular and truncated $A$.}
The spectral theory requires a square surrogate pair, but the original
operator used for transfer may be rectangular. In Experiment~3,
selection uses the mirror-extended data
$b_{\mathrm{sel}}=\mathsf{Mir}\,b_{\mathrm{orig}}$,
where $b_{\mathrm{orig}}$ is the recorded window. The surrogate
$\widetilde A$, perturbation $S$ and penalty $L_{\mathrm{sel}}=I$
act on the extended grid (\cref{sec:real-data-setup}). The discrepancy
is normalized by $\|b_{\mathrm{sel}}\|_2$, and
\cref{sec:spectral-theory,sec:controller} apply to this square pair.

Transfer uses the recorded data and the original operator, whose unknown
lies on the field of view extended by the PSF radius:
\(
x_{\lambda_{\mathrm{SSDM}}}(A)
=
(A^*A+\lambda_{\mathrm{SSDM}}I)^{-1}A^*b_{\mathrm{orig}}.
\)
Thus selection and transfer share the scalar parameter and use identity
penalties on their respective grids, but their data and unknowns have
different dimensions. In this setting, $\widetilde A$ is not an
approximation of $A$ in common coordinates, so
$A-\widetilde A$ is not directly defined; operator mismatch is assessed
on a common restricted domain instead (\cref{app:exp3-setup}).

% \paragraph{Rectangular and truncated $A$.}
% The restriction to square $A$ in Step~1 serves the spectral theory, which
% concerns only the surrogate pair, and that pair is square by construction.
% In Experiment~3, the original operator is rectangular and acts on grids
% different from those of the surrogate pair. Write $b_{\mathrm{orig}}$ for the
% recorded data (the native window) and $b_{\mathrm{sel}}=\mathsf{Mir}\,b_{\mathrm{orig}}$
% for its \emph{mirror} extension, the data of the selection pair. The surrogate
% $\widetilde A$, the perturbation $S$, and the selection penalty
% $L_{\mathrm{sel}}=I$ act on the grid of $b_{\mathrm{sel}}$
% (\cref{sec:real-data-setup}). Selection solves $D_b(\lambda)=\tau$ with
% $b_{\mathrm{sel}}$ and $\|b_{\mathrm{sel}}\|_2$, and
% \cref{sec:spectral-theory,sec:controller} apply to this square pair
% unchanged. Transfer uses the recorded data and the unknown of $A$, which
% lives on the field of view extended by the PSF radius:
% $x_{\lambda_{\mathrm{SSDM}}}(A)=(A^*A+\lambda_{\mathrm{SSDM}}L_{\mathrm{orig}}^*L_{\mathrm{orig}})^\dagger
% A^*b_{\mathrm{orig}}$ with $L_{\mathrm{orig}}=I$ on that grid. In this setting
% ``the same data'' and ``the same penalty'' in Steps~2 and~4 mean the same up to
% the extension $\mathsf{Mir}$ and the identity on the respective grid; they are
% not matrices of equal dimension, and only the scalar $\lambda$ is shared. What
% is lost is the interpretation of $\widetilde A$ as an approximation of $A$ in
% the same coordinates, so the mismatch cannot be evaluated directly on the
% full field (\cref{app:exp3-setup}).

\section{The data-independent structure of the discrepancy}
\label{sec:spectral-theory}

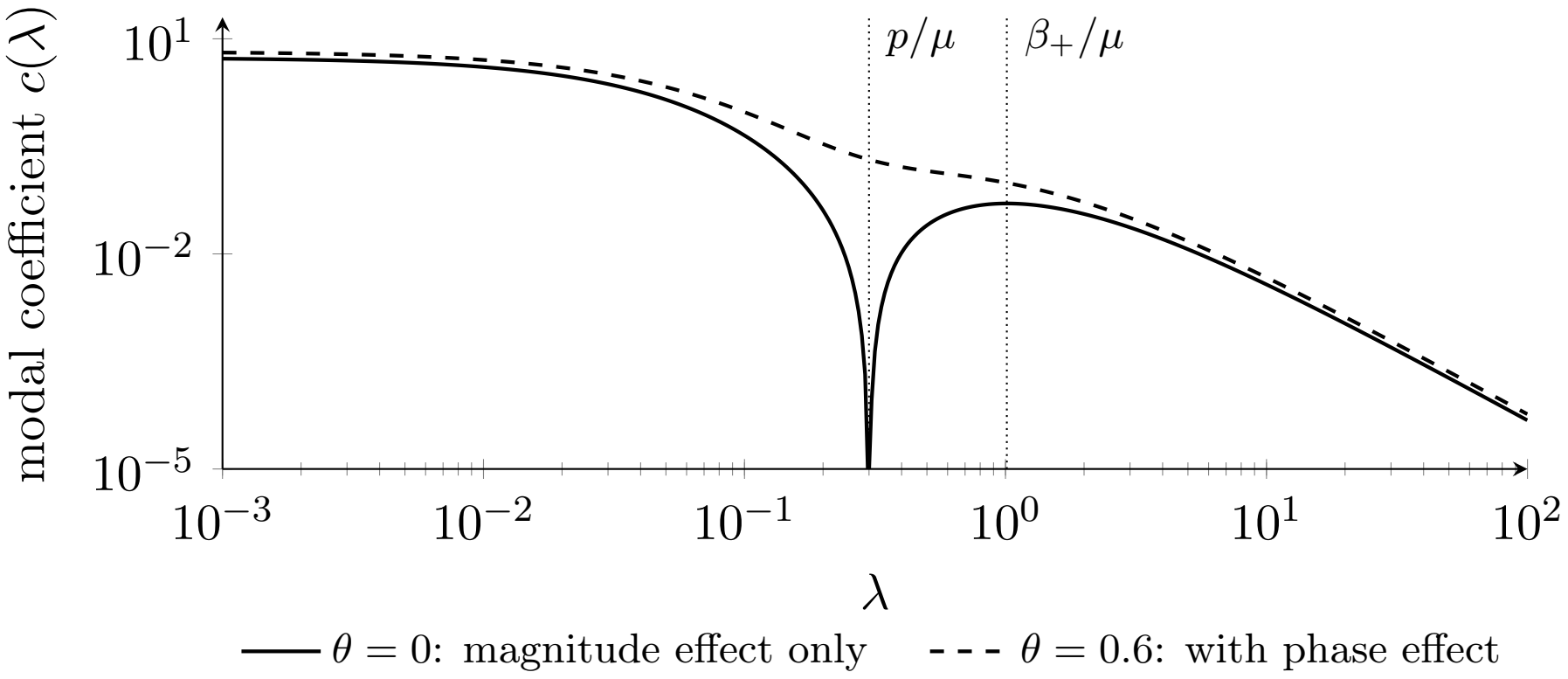
\begin{figure}[tb]
\centering
\begin{tikzpicture}
\begin{axis}[
  width=0.72\textwidth, height=5cm,
  xmode=log, ymode=log,
  xmin=1e-3, xmax=1e2, ymin=1e-5, ymax=20,
  xlabel={$\lambda$}, ylabel={modal coefficient $c(\lambda)$},
  axis lines=left, legend cell align=left, legend columns=2,
  legend style={font=\footnotesize,draw=none,fill=none,
    at={(0.5,-0.32)},anchor=north,/tikz/every even column/.append style={column sep=1em}},
  every axis plot/.append style={thick}]
  \addplot[domain=1e-3:1e2, samples=401, black]
    {(1/(1+x)-0.3/(0.09+x))^2};
  \addlegendentry{$\theta=0$: magnitude effect only}
  \addplot[domain=1e-3:1e2, samples=401, black, dashed]
    {(1/(1+x)-0.3/(0.09+x))^2+2*(1-cos(deg(0.6)))*0.3/((1+x)*(0.09+x))};
  \addlegendentry{$\theta=0.6$: with phase effect}
  \draw[densely dotted] (axis cs:0.3,1e-5) -- (axis cs:0.3,20);
  \draw[densely dotted] (axis cs:1.012,1e-5) -- (axis cs:1.012,20);
  \node[font=\footnotesize,anchor=south west] at (axis cs:0.3,3) {$p/\mu$};
  \node[font=\footnotesize,anchor=south west] at (axis cs:1.012,3) {$\beta_+/\mu$};
\end{axis}
\end{tikzpicture}
\caption{The modal coefficient $c$ of \eqref{eq:general-ck} for $a=1$,
$s=0.3$, $\mu=1$. Without relative phase (solid), $c=h^2$ dips to zero at
$\lambda=p/\mu$, rises, and decreases only beyond $\beta_+/\mu$
(\cref{lem:scalar-filter-difference}). With relative phase $\theta=0.6$
(dashed), the nonincreasing phase term of \cref{lem:phase-decomposition}
fills the dip; here condition \eqref{eq:phase-dominance} holds, and $c$
decreases on all of $(0,\infty)$.}
\label{fig:modal-geometry}
\end{figure}

Because both models are diagonal in the same basis, the squared discrepancy
is a weighted sum of scalar functions, one per mode, and the weights are the
data energies in that basis. Everything in this section follows from the
behavior of a single mode, where two effects compete. The two models
attenuate the mode by different amounts. This magnitude effect makes the
modal discrepancy dip, then rise, then decay. The two models may also shift
the mode's phase differently, and this phase effect contributes a term that
decreases monotonically (\cref{fig:modal-geometry}). The main theorem
(\cref{thm:main-spectral-general}) uses the magnitude effect to obtain the
global shape of the discrepancy and a closed-form monotonicity threshold.
Keeping the phase effect gives a sharper threshold, which cannot be improved
without looking at the data (\cref{thm:phase-aware-threshold}).

\subsection{Assumptions and spectral representation}
\label{subsec:general-phase-formula}
\label{subsec:tractability-conditions}

Let $\widetilde A$, $S$, and $L$ have the common unitary diagonalization
\eqref{eq:simultaneous-diagonalization}, so that
$\widetilde A_S=Q\diag(\nu\odot\gamma)Q^*$, where $\odot$ denotes element-wise or Hadamard product. Set
\[
\widehat b:=Q^*b,\qquad
\mathcal K:=\{k:\widehat b_k\ne0\},\qquad
\mu_k:=|\nu^L_k|^2 .
\]
We call the indices in $\mathcal K$ \emph{active modes}; only these enter
the reconstructions. Write
\[
\nu_k=a_ke^{i\varphi_k},\qquad
\gamma_k=\rho_ke^{i\theta_k},\qquad
s_k:=|\nu_k\gamma_k|=a_k\rho_k,
\]
so that $a_k$ and $s_k$ are the modal magnitudes of the two models and
$\theta_k$ is their relative phase. The phase of a zero number may be
chosen arbitrarily: every phase-dependent term below is multiplied by the
corresponding magnitude and vanishes.

\begin{definition}[Mode-wise contraction]
\label{def:modewise-contraction}
The perturbation $S$ is a \emph{mode-wise contraction} if $|\gamma_k|\le1$
for every $k$. This imposes no phase restriction and implies
$0\le s_k\le a_k$.
\label{eq:contraction-def}
\label{eq:nonnegative-multiplier-assumption}
\end{definition}

\begin{assumption}[Standing assumptions]
\label{ass:standing}
$S$ is a mode-wise contraction, and $\forall k\in\mathcal K, \mu_k>0$.
\end{assumption}

Contraction expresses that the second model does not amplify any mode. It is
the motivating case, and the half-pixel pair satisfies it, but the results
below do not need it. The discrepancy, its phase term and the thresholds are
all symmetric in the two modal magnitudes. For any two models diagonalized by the same
unitary transform, with the same positive active penalties, the statements
therefore hold after ordering each pair of magnitudes and replacing
$\mathcal J$ by $\{k\in\mathcal K:a_ks_k>0,\ a_k\ne s_k\}$
(\cref{rem:beyond-contraction}). We keep the notation $s_k\le a_k$ throughout.
Positivity of the penalty is needed only on
active modes. Under these assumptions the normal matrices of both models are
diagonalized by $Q$, their entries are at least $\lambda\mu_k>0$ on active
modes, and inactive modes contribute nothing to either reconstruction.
Unitarity of $Q$ therefore gives the following representation.

\begin{proposition}[Spectral representation]
\label{prop:general-phase-formula}
For $\lambda>0$,
\begin{equation}
\Phi(\lambda)=\sum_{k\in\mathcal K}c_k(\lambda)|\widehat b_k|^2,\qquad
c_k=\phi_k^2+(\phi_k^S)^2-2\phi_k\phi_k^S\cos\theta_k,
\label{eq:general-ck}
\end{equation}
where
\begin{equation}
\phi_k(\lambda):=\frac{a_k}{a_k^2+\lambda\mu_k},\qquad
\phi_k^S(\lambda):=\frac{s_k}{s_k^2+\lambda\mu_k}.
\label{eq:phi-def}
\end{equation}
\end{proposition}
\begin{proof}
The $k$th coefficient of $Q^*(\widetilde x_\lambda-\widetilde x_\lambda^S)$
is $\bigl(\overline{\nu_k}/(a_k^2+\lambda\mu_k)
-\overline{\nu_k\gamma_k}/(s_k^2+\lambda\mu_k)\bigr)\widehat b_k
=e^{-i\varphi_k}(\phi_k-\phi_k^Se^{-i\theta_k})\widehat b_k$, whose squared
modulus is $c_k|\widehat b_k|^2$.
\end{proof}

The filters $\phi_k$ and $\phi_k^S$ are the regularized inverses of the two
modal magnitudes. For $t\ge0$ write $t^\dagger:=1/t$ if $t>0$ and
$t^\dagger:=0$ if $t=0$. On active modes we set $\phi_k(0):=a_k^\dagger$ and
$\phi_k^S(0):=s_k^\dagger$; since $\mu_k>0$, these are the limits as
$\lambda\downarrow0$, so the extended functions are continuous on
$[0,\infty)$.
For the normalized discrepancy, $\Phi_b=\sum_kc_kw_k$ with weights
$w_k:=|\widehat b_k|^2/\|b\|_2^2$ that sum to one. Every result below holds
verbatim for $\Phi_b$ and $D_b$ after replacing $|\widehat b_k|^2$ by $w_k$.
\label{cor:data-normalization}%
\label{eq:normalized-weights}%

\subsection{One mode: magnitude and phase}
\label{subsec:scalar-filter-geometry}
\label{subsec:phase-domination}

Fix an active mode and suppress its index. For $a\ge s\ge0$ and $\mu>0$ let
\begin{equation}
h(\lambda):=\frac{a}{a^2+\lambda\mu}-\frac{s}{s^2+\lambda\mu},
\label{eq:difference-filter}
\end{equation}
the difference of the two filters, which carries the magnitude effect. Under
phase alignment ($\theta=0$) the modal coefficient is $c=h^2$. In the proofs
we use $\beta:=\lambda\mu$, $p:=as$, and
\begin{equation}
N(\beta):=(a^2+\beta)(s^2+\beta),\qquad
P_0(\beta):=\beta^2-2p\beta-p(a^2+p+s^2)=(\beta-\beta_+)(\beta-\beta_-),
\label{eq:N-P0-def}
\end{equation}
\begin{equation}
\beta_\pm:=p\pm(a+s)\sqrt p .
\label{eq:beta-pm}
\end{equation}
If $p>0$, then $\beta_-<0<p<\beta_+$.

\begin{lemma}[Magnitude effect]
\label{lem:scalar-filter-difference}
If $a=s$, then $h\equiv0$. If $a>s=0$, then $h^2$ is strictly decreasing on
$(0,\infty)$. If $a>s>0$, then $h(\lambda)^2<h(0)^2=(1/a-1/s)^2$ for every
$\lambda>0$, and $h^2$ is strictly decreasing on $(0,p/\mu)$, strictly
increasing on $(p/\mu,\beta_+/\mu)$, and strictly decreasing on
$(\beta_+/\mu,\infty)$. For all $a\ge s\ge0$,
\begin{equation}
h(\lambda)^2=\frac{(a-s)^2}{\lambda^2\mu^2}+O(\lambda^{-3}),
\qquad\lambda\to\infty.
\label{eq:scalar-filter-asymptotics}
\end{equation}
\end{lemma}
\begin{proof}
The cases $a=s$ and $a>s=0$ are immediate. For $a\ge s\ge0$ and $\beta>0$, a
common denominator gives
\begin{equation}
\widetilde h(\beta):=h(\beta/\mu)=\frac{(a-s)(\beta-p)}{N(\beta)},\qquad
\widetilde h'(\beta)=-\frac{(a-s)P_0(\beta)}{N(\beta)^2},
\label{eq:htilde}
\end{equation}
where the derivative uses $N'(\beta)=a^2+s^2+2\beta$ and the identity
$N(\beta)-(\beta-p)(a^2+s^2+2\beta)=-P_0(\beta)$.
\label{eq:htilde-derivative}%
Now let $a>s>0$. The bound $|\widetilde h(\beta)|<|\widetilde h(0)|$ is
equivalent to $p|\beta-p|<N(\beta)$, and $N(\beta)-p|\beta-p|$ equals
$\beta(a^2+p+s^2+\beta)$ for $0<\beta<p$ and
$2p^2+\beta(a^2-p+s^2)+\beta^2$ for $\beta\ge p$, both positive. By
\eqref{eq:htilde}, $\widetilde h<0$ on $(0,p)$ and $\widetilde h>0$ on
$(p,\infty)$, while $\widetilde h'>0$ on $(0,\beta_+)$ and $\widetilde h'<0$
on $(\beta_+,\infty)$. Since $(\widetilde h^2)'=2\widetilde h\widetilde h'$
and $\beta=\lambda\mu$, this gives the three monotonicity intervals.
Finally, $a/(a^2+\lambda\mu)=a/(\lambda\mu)+O(\lambda^{-2})$, also for
$a=0$, which gives \eqref{eq:scalar-filter-asymptotics}.
\end{proof}

\begin{lemma}[Phase effect]
\label{lem:phase-decomposition}
For $a\ge s\ge0$, $\mu>0$, and $\theta\in\R$, let $\phi=a/(a^2+\lambda\mu)$,
$\phi^S=s/(s^2+\lambda\mu)$, and $c=\phi^2+(\phi^S)^2-2\phi\phi^S\cos\theta$.
Then
\begin{equation}
c(\lambda)=h(\lambda)^2+2(1-\cos\theta)\phi(\lambda)\phi^S(\lambda),
\qquad\lambda\ge0.
\label{eq:phase-decomposition}
\end{equation}
The phase term $2(1-\cos\theta)\phi\phi^S$ is nonnegative and
nonincreasing; it is strictly decreasing if $as>0$ and
$\theta\not\equiv0\pmod{2\pi}$, and vanishes identically otherwise. In
particular $c\ge h^2$ and
\begin{equation}
c'(\lambda)\le(h^2)'(\lambda),\qquad\lambda>0.
\label{eq:phase-derivative-bound}
\end{equation}
\end{lemma}
\begin{proof}
Adding and subtracting $2\phi\phi^S$ gives \eqref{eq:phase-decomposition}.
Both filters are nonnegative with derivatives
$\phi'=-\mu a/(a^2+\lambda\mu)^2\le0$ and
$(\phi^S)'=-\mu s/(s^2+\lambda\mu)^2\le0$, so
$c'-(h^2)'=2(1-\cos\theta)(\phi'\phi^S+\phi(\phi^S)')\le0$, strictly if
$as>0$ and $\cos\theta<1$. If $as=0$ or $\cos\theta=1$, the term vanishes.
\end{proof}

\subsection{The global shape of the discrepancy}
\label{subsec:main-theorem}

\begin{theorem}[Shape of the SSDM discrepancy]
\label{thm:main-spectral-general}
\label{prop:boundary-maximum}
\label{thm:asymptotic-decay}
\label{cor:asymptotic-rate}
Let \cref{ass:standing} hold, and suppose some active mode $k_0$ satisfies
\begin{equation}
\nu_{k_0}(1-\gamma_{k_0})\ne0,
\label{eq:active-model-nondegeneracy}
\end{equation}
that is, the perturbation changes the model on at least one mode excited by
the data. Let $\mathcal J:=\{k\in\mathcal K:a_k>s_k>0\}$,
$\beta_{k,+}:=a_ks_k+(a_k+s_k)\sqrt{a_ks_k}$, and
\[
\lambda^\ddagger:=\max_{k\in\mathcal J}\beta_{k,+}/\mu_k
\qquad(\lambda^\ddagger:=0\text{ if }\mathcal J=\emptyset).
\]
Then $\Phi\in C^\infty((0,\infty))$, and:
\begin{enumerate}[label=(\roman*)]
\item \emph{Boundary maximum.} The continuous extension to $\lambda=0$ is
$\Phi(0)=\sum_{k\in\mathcal K}\bigl[(a_k^\dagger)^2+(s_k^\dagger)^2
-2a_k^\dagger s_k^\dagger\cos\theta_k\bigr]|\widehat b_k|^2$, and
$\Phi(\lambda)<\Phi(0)$ for every $\lambda>0$.
\item \emph{Decreasing tail.} $\Phi'(\lambda)<0$ for every
$\lambda>\lambda^\ddagger$. There $D>0$ is smooth and $D'=\Phi'/(2D)<0$.
\item \emph{Decay.} As $\lambda\to\infty$,
\begin{equation}
\begin{gathered}
\Phi(\lambda)=C\lambda^{-2}+O(\lambda^{-3}),\qquad
D(\lambda)=\sqrt C\,\lambda^{-1}+O(\lambda^{-2}),\\
C:=\sum_{k\in\mathcal K}\frac{|\nu_k(1-\gamma_k)|^2}{\mu_k^2}|\widehat b_k|^2>0.
\end{gathered}
\label{eq:asymptotic-numerator-expansion}
\end{equation}
\end{enumerate}
If \eqref{eq:active-model-nondegeneracy} fails at every active mode, then
$\Phi\equiv0$.
\end{theorem}

Parts (i) and (iii) say that every target $0<\tau<D(0)$ is attained
somewhere on $(0,\infty)$. They do not say where, or how often: between
the boundary and the tail, the dip and rise of
\cref{lem:scalar-filter-difference} can survive in the sum. The value
$\Phi(0)$ is dominated by weakly observed modes and says nothing about the
range used for selection. The threshold $\lambda^\ddagger$ is the point
beyond which the magnitude effect of every active mode has stopped rising.

\begin{proof}
\emph{Boundary value.} Each $c_k$ is smooth on $(0,\infty)$ and continuous
at $0$, which gives the formula for $\Phi(0)$. For an active mode with
$a>s>0$, \cref{lem:scalar-filter-difference} gives $h(\lambda)^2<h(0)^2$
and \cref{lem:phase-decomposition} a nonincreasing phase term, so
$c(\lambda)<c(0)$. If $a>s=0$, then $c=\phi^2$ is strictly decreasing. If
$a=s>0$, then $h\equiv0$ and
\begin{equation}
c(\lambda)=2(1-\cos\theta)\frac{a^2}{(a^2+\lambda\mu)^2},
\label{eq:pure-phase-mode}
\end{equation}
which is strictly decreasing unless $\gamma=1$, when $c\equiv0$. If $a=0$,
then $c\equiv0$. Thus $c_k(\lambda)\le c_k(0)$ on every active mode, strictly
when $\nu_k(1-\gamma_k)\ne0$, and the positive weight of $k_0$ gives (i).

\emph{Tail.} For $k\in\mathcal J$ and $\lambda>\beta_{k,+}/\mu_k$, the two
lemmas give $c_k'\le(h_k^2)'<0$. The case list above shows that every other
active mode has $c_k'\le0$ on $(0,\infty)$, with $c_{k_0}'<0$. Hence
$\Phi'<0$ for $\lambda>\lambda^\ddagger$. Since $\Phi\ge0$ and $\Phi$
strictly decreases there, $\Phi>0$ there, and $D=\sqrt\Phi$ is smooth.

\emph{Decay.} Expanding
$\phi_k=a_k/(\lambda\mu_k)-a_k^3/(\lambda\mu_k)^2+O(\lambda^{-3})$ and
similarly $\phi_k^S$ in \eqref{eq:general-ck} gives
$c_k=(a_k^2+s_k^2-2a_ks_k\cos\theta_k)/(\lambda\mu_k)^2+O(\lambda^{-3})$,
and $a_k^2+s_k^2-2a_ks_k\cos\theta_k=|\nu_k(1-\gamma_k)|^2$. Summing gives
(iii), with $C>0$ because of $k_0$. If $\nu_k(1-\gamma_k)=0$ for every
active $k$, the two filters coincide on every active mode and $\Phi\equiv0$.
\end{proof}

\subsection{A sharper threshold, and why it cannot be improved without the data}
\label{subsec:phase-aware-threshold}

The threshold $\lambda^\ddagger$ ignores the phase term, which pulls every
mode downward. For perturbations that act mainly as phase shifts, such as
the half-pixel shift, this is wasteful. Keeping it requires the exact modal
derivative. For an active mode let $t_k:=1-\cos\theta_k\in[0,2]$,
$d_k:=(a_k-s_k)^2$, $p_k:=a_ks_k$, let $N_k$ and $P_{0,k}$ be
\eqref{eq:N-P0-def} with $(a,s)=(a_k,s_k)$, and define the cubic
\begin{equation}
P_k(\beta):=d_k(\beta-p_k)P_{0,k}(\beta)
+t_k\,p_k\,(a_k^2+s_k^2+2\beta)N_k(\beta).
\label{eq:Pk-def}
\end{equation}

\begin{lemma}[Exact modal derivative]
\label{lem:exact-modal-derivative}
For every active mode and every $\lambda>0$, with $\beta=\lambda\mu_k$,
\begin{equation}
c_k'(\lambda)=-\frac{2\mu_k\,P_k(\beta)}{N_k(\beta)^3}.
\label{eq:exact-modal-derivative}
\end{equation}
$P_k\equiv0$ if and only if $\nu_k(1-\gamma_k)=0$; otherwise its leading
coefficient is $|\nu_k(1-\gamma_k)|^2>0$.
\end{lemma}
\begin{proof}
By \eqref{eq:phase-decomposition} and \eqref{eq:htilde},
$c=\widetilde h^2+2tp/N$ as a function of $\beta$. Since
$(\widetilde h^2)'=-2d(\beta-p)P_0/N^3$ and $(p/N)'=-p(a^2+s^2+2\beta)/N^2$,
$dc/d\beta=-2P/N^3$, and $d\beta/d\lambda=\mu_k$. The leading coefficient of
$P$ is $d+2tp=|a-se^{i\theta}|^2=|\nu(1-\gamma)|^2$; if it vanishes, then
$d=0$ and $tp=0$, so $P\equiv0$. The monomial coefficients are listed in
\cref{app:certificate-computation}.
\end{proof}

A mode increases exactly where its cubic is negative. Its largest such
point is
\begin{equation}
\beta_k^\theta:=\sup\{\beta>0:P_k(\beta)\le0\}\quad(\sup\emptyset:=0;\
\beta_k^\theta:=0\text{ if }P_k\equiv0),
\qquad
\lambda^\theta:=\max_{k\in\mathcal K}\beta_k^\theta/\mu_k ,
\label{eq:beta-theta-def}
\end{equation}
and $\beta_k^\theta$ is the largest positive real zero of $P_k$, or zero.
We say that mode $k$ satisfies the \emph{sign-change condition} if
$P_k$ actually changes sign at $\beta_k^\theta$, that is,
\begin{equation}
\beta_k^\theta=\sup\{\beta>0:P_k(\beta)<0\}.
\label{eq:modal-sign-change}
\end{equation}
This excludes only a terminal zero at which $P_k$ touches zero without
becoming negative. It holds automatically for phase-aligned modes
($\theta_k\equiv0$), which set the threshold for the half-pixel pair on
square grids (\cref{cor:half-pixel-thresholds}). On other grids the
maximizing mode can be nearly, but not exactly, phase-aligned. There the condition is not automatic and must be checked (~\cref{app:rect-grid}).

\begin{theorem}[Phase-aware threshold; its optimality]
\label{thm:phase-aware-threshold}
Assume the hypotheses of \cref{thm:main-spectral-general}.
\begin{enumerate}[label=(\roman*)]
\item \emph{Decreasing tail.} $\Phi'(\lambda)<0$ for
$\forall\lambda>\lambda^\theta$. \cref{thm:main-spectral-general}(ii) conclusions hold on $(\lambda^\theta,\infty)$.
\item \emph{Comparison.} $\beta_k^\theta=0$ if $s_k=0$ or $a_k=s_k$. If
$a_k>s_k>0$, then $\beta_k^\theta\le\beta_{k,+}$, with equality if and only
if $\theta_k\equiv0$, and $\beta_k^\theta$ is nonincreasing in $t_k$.
Hence $\lambda^\theta\le\lambda^\ddagger$.
\item \emph{Phase dominance.} If $a_k>s_k>0$ and
\begin{equation}
16\,a_ks_k(1-\cos\theta_k)\ge(a_k-s_k)^2,
\label{eq:phase-dominance}
\end{equation}
then $\beta_k^\theta=0$: the mode decreases strictly on $(0,\infty)$. The
constant $16$ cannot be increased.
\item \emph{Optimality.} If $P_k(\beta_1)<0$ for some active $k$ and
$\beta_1>0$, then for any fixed positive weights on the other active modes
and all sufficiently large weight $\omega_k$ on mode $k$, the weighted sum
$\Phi_\omega=\sum_jc_j\omega_j$ is increasing at $\lambda_1=\beta_1/\mu_k$.
Consequently, if a mode attaining the maximum in \eqref{eq:beta-theta-def}
satisfies the sign-change condition, no $\bar\lambda<\lambda^\theta$ has the
property that $\Phi_\omega$ decreases on $(\bar\lambda,\infty)$ for every
positive weight vector $\omega$ on $\mathcal K$.
\end{enumerate}
\end{theorem}
\begin{proof}[Proof of (i)]
For $\lambda>\lambda^\theta$ and active $k$, $\beta=\lambda\mu_k>\beta_k^\theta$,
so $P_k(\beta)>0$ unless $P_k\equiv0$, and
\cref{lem:exact-modal-derivative} gives $c_k'(\lambda)\le0$, strictly for
$k_0$. Parts (ii)--(iv) are proved in \cref{app:phase-aware-proofs}.
\end{proof}

Part (iv) concerns arbitrary positive modal weights. Real data constrain
the weights by conjugate symmetry, but for a real Fourier-diagonal pair the
modes $k$ and $-k$ carry the same $(a_k,s_k,t_k,\mu_k)$. Placing equal
weight on both realizes the construction with real data. The statement is
therefore a limit on what any argument that ignores the data weights
can achieve.

\subsection{What the data-independent theory does not reach}
\label{rem:two-levels}
\label{rem:scope-guarantees}

The results of this section concern the surrogate pair and hold for every
data vector satisfying the active-mode hypotheses. They do not require $b$
to be generated by $\widetilde A$, and they say nothing about the accuracy of
the transferred reconstruction. The results explain the global shape of the
discrepancy, but their optimality is relative to thresholds that
ignore the data weights, and that class could be outside the range of
interest. For the half-pixel deblurring pair the threshold $\lambda^\theta$
lies between about $1.6$ and $3$ (\cref{cor:half-pixel-thresholds},
\cref{app:rect-grid}), whereas the parameters selected in our
experiments lie between $3\times10^{-4}$ and $0.12$. By
\cref{thm:phase-aware-threshold}(iv), no data-independent argument can close
this gap. Monotonicity in the selection range is a property of the data, which the next section certifies. Values of the discrepancy or
its derivative at finitely many sample points cannot do this, since they do
not exclude extrema between the samples.

\section{Certified target matching}
\label{sec:controller}

Step~3 of SSDM asks for $\lambda$ with $D_b(\lambda)=\tau$ . By
\cref{thm:phase-aware-threshold}(i), every target
$0<\tau<D_b(\lambda^\theta)$ has a unique root in
$(\lambda^\theta,\infty)$.
\label{prop:controller-tail-interval}%
Below $\lambda^\theta$, where selection happens in practice, the sign of
$\Phi_b'$ depends on the data, and a target may be met several times
(\cref{app:controller-multiple-crossings}). This section certifies strict
decrease on a compact interval by a finite computation and then locates
the root. We formulate the derivative bounds for $\Phi$ and obtain their
normalized counterparts by division by $\|b\|_2^2$.

\subsection{A finite certificate of strict decrease}
\label{subsec:controller-intervals-coordinates}
\label{subsec:certified-intervals}

By \cref{lem:exact-modal-derivative}, $c_k'(\lambda)$ is a positive
multiple of $-P_k(\lambda\mu_k)$. On a short interval, an upper bound for
$c_k'$ therefore follows from a lower bound for the cubic $P_k$ and the
monotonicity of $N_k$. Fix a partition
$0<\bar\lambda_0<\bar\lambda_1<\cdots<\bar\lambda_J$ and, for each cell $j$
and active mode $k$, let
\[
\beta_{k,j}^-:=\mu_k\bar\lambda_{j-1},\qquad
\beta_{k,j}^+:=\mu_k\bar\lambda_j,\qquad
m_{k,j}:=\min_{\beta\in[\beta_{k,j}^-,\beta_{k,j}^+]}P_k(\beta),
\]
and define
\begin{equation}
u_{k,j}:=\begin{cases}
-2\mu_k\,m_{k,j}/N_k(\beta_{k,j}^+)^3,&m_{k,j}\ge0,\\
-2\mu_k\,m_{k,j}/N_k(\beta_{k,j}^-)^3,&m_{k,j}<0,
\end{cases}
\qquad
U_j:=\sum_{k\in\mathcal K}u_{k,j}|\widehat b_k|^2 .
\label{eq:certificate-upper}
\end{equation}
For the normalized discrepancy, put
\[
U_{b,j}:=\frac{U_j}{\|b\|_2^2}
=\sum_{k\in\mathcal K}u_{k,j}w_k,
\qquad
w_k:=\frac{|\widehat b_k|^2}{\|b\|_2^2}.
\]
Since $b\ne0$, $U_{b,j}$ and $U_j$ have the same sign.

\begin{proposition}[Data-dependent interval certificate]
\label{prop:interval-certificate}
Under \cref{ass:standing}, with $b\ne0$, the bounds
$\Phi'(\lambda)\le U_j$ and $\Phi_b'(\lambda)\le U_{b,j}$ hold for every
$\lambda\in[\bar\lambda_{j-1},\bar\lambda_j]$. If $U_j<0$ for all $j$,
then on $I:=[\bar\lambda_0,\bar\lambda_J]$ both $D$ and $D_b$ are
positive, smooth, and strictly decreasing. Every target
\[
D_b(\bar\lambda_J)\le\tau\le D_b(\bar\lambda_0)
\]
then has exactly one root of $D_b(\lambda)=\tau$ in $I$.
If moreover $\bar\lambda_J>\lambda^\theta$, both discrepancies decrease
strictly on $[\bar\lambda_0,\infty)$. Conversely, if $\Phi'<0$ on $I$,
equivalently $\Phi_b'<0$ there, then $U_j<0$ for all $j$ once the mesh is
fine enough.
\end{proposition}

The proof is given in \cref{app:certificate-proof}. Dividing by $\|b\|_2^2$ and $\|b\|_2$ gives us the normalized statements. The converse matters in practice: the
certificate succeeds whenever the derivative is strictly negative
throughout the compact interval and the partition is sufficiently fine.

% \begin{remark}[Computing the certificate]
% \label{rem:certificate-computation}
% The extrema of a cubic on an interval lie at the endpoints or at the real
% zeros of its derivative, so each $m_{k,j}$ is available in closed form.
% Once the spectra and data coefficients are known, the certificate costs
% $O(Jn)$ operations. Modes with $P_k>0$ on $(0,\infty)$ contribute
% strictly negative terms to every upper sum. They may therefore be omitted
% from the sign test: omission increases the upper bound and preserves
% sufficiency. The refinement guarantee in
% \cref{prop:interval-certificate}, however, applies to the full sum;
% omitting decreasing modes may prevent the test from succeeding.

% A computer-assisted proof would require guaranteed enclosures of every
% quantity involved, including the polynomial extrema, data weights, and
% sums, for example through outward-rounded interval arithmetic. Our
% computations evaluate the bounds in ordinary floating point. They are
% therefore numerical verifications of the hypotheses, not
% computer-assisted proofs, and we call them \emph{numerically verified}
% throughout.
% \end{remark}

\subsection{Finding the root}
\label{subsec:bracketed-root}
\label{prop:bracketed-root}

On a certified interval whose endpoint values bracket the target,
continuity and strict decrease give a unique root. A method that
maintains a sign-change bracket for $D_b-\tau$ keeps this root in its
bracket. Convergence follows when the bracket lengths tend to zero. An
endpoint that attains the target is the root. Such methods need only the sign of $D_b-\tau$. Once the spectra and the data coefficients $\widehat b=Q^*b$ are available, each evaluation of $D_b$ is a spectral sum of $O(n)$ operations that needs no further transform, and working in $\log\lambda$ treats parameters spread over many decades evenly~(\Cref{alg:ssdm}). A failed test is inconclusive: it does not show that the discrepancy is nonmonotone. If the derivative is strictly negative throughout the interval, sufficiently fine refinement will succeed by \cref{prop:interval-certificate}.

\begin{algorithm}[tb]
\caption{SSDM parameter selection with an interval certificate}
\label{alg:ssdm}
\begin{algorithmic}[1]
\Require data $b\neq0$; spectra $\nu_k,\gamma_k,\mu_k$ in the common
basis $Q$, satisfying \cref{ass:standing}; target $\tau>0$; search
interval $[\bar\lambda_0,\bar\lambda_J]$ with a log partition;
finite refinement budget $R$
\State compute $\widehat b=Q^*b$,
$\mathcal K=\{k:\widehat b_k\ne0\}$, and
$w_k=|\widehat b_k|^2/\|b\|_2^2$
\State compute $U_{b,j}$ using \eqref{eq:certificate-upper} with $w_k$
in place of $|\widehat b_k|^2$, for $j=1,\ldots,J$
\State set $r=0$
\While{some $U_{b,j}\ge0$}
  \If{$r=R$}
    \State \Return \textsc{inconclusive}: no interval certificate obtained
  \EndIf
  \State refine the partition near the failing cells
  \State recompute the affected bounds and set $r=r+1$
\EndWhile
\If{$\tau\notin[D_b(\bar\lambda_J),D_b(\bar\lambda_0)]$}
  \State \Return \textsc{unattained}: target outside the interval range
\EndIf
\State locate the unique root $\lambda_{\mathrm{SSDM}}$ of
$D_b(\lambda)=\tau$ to the prescribed tolerance by bisection or a
safeguarded bracketing method in $\log\lambda$, returning an endpoint
if it attains the target
\State \Return $x_{\lambda_{\mathrm{SSDM}}}(A)$, e.g.\ by conjugate
gradients on the normal equations, preconditioned with the surrogate
\end{algorithmic}
\end{algorithm}

\section{Half-pixel deblurring}
\label{sec:problem-setup}

In our running example, the surrogate is a periodic
(BCCB) blur, and the second model shifts the image by half a pixel along
the diagonal. Two questions from the general theory become concrete. Where
does the data-independent threshold sit? Near $3$ on square grids, for a
reason that can be read off the geometry of the shift
(\cref{cor:half-pixel-thresholds}). And what does the discrepancy measure?
A directional roughness of the reconstruction (\cref{subsec:roughness}).

Let $\widetilde A=p*$ be periodic convolution with a kernel $p$ on an
$n_1\times n_2$ grid, $n=n_1n_2$, with images identified with their
column-vectorizations, and let $L=I_n$, so $\mu_k=1$. Let $F$ be the
unitary Fourier synthesis matrix, so $Q=F$ and $\widehat b=F^*b$. With
periodic indices, $(p*x)_{i,j}=\sum_{r,t}p_{r,t}x_{i-r,j-t}$ and
\[
\nu_k=\sum_{r,t}p_{r,t}e^{-i(r\omega_{1,k}+t\omega_{2,k})}=\sqrt n\,(F^*p)_k .
\]
Convolution eigenvalues thus use the unnormalized DFT of the kernel indexed
from its origin, while data coefficients use the unitary DFT
(\cref{app:fourier-implementation}). The bilinear half-pixel operator for
the shift $v=(h_{\mathrm{px}}/2,h_{\mathrm{px}}/2)$, with pixel side
$h_{\mathrm{px}}$, is
\[
(Sx)_{i,j}=\tfrac14(x_{i,j}+x_{i-1,j}+x_{i,j-1}+x_{i-1,j-1})=(d_v*x)_{i,j},
\]
where $d_v$ has weights $1/4$ at offsets $(0,0),(1,0),(0,1),(1,1)$. It
evaluates the bilinear interpolant at $(i-\frac12,j-\frac12)$, representing
$x(\cdot-v)$, and is generally not a unitary translation. The second model
is $\widetilde A_S=\widetilde AS=(p*d_v)*$.
\label{eq:perturbed-forward-model}

\begin{proposition}[Bilinear half-pixel multiplier]
\label{prop:half-pixel-multiplier}
Given indices $k_j\in\{-\lfloor n_j/2\rfloor,\ldots,\lceil n_j/2\rceil-1\}$
and frequencies $\omega_{j,k}=2\pi k_j/n_j\in[-\pi,\pi)$, 
$S=F\diag(\gamma)F^*$ with
\begin{equation}
\begin{gathered}
\gamma_k=\tfrac14(1+e^{-i\omega_{1,k}})(1+e^{-i\omega_{2,k}})
=\rho_ke^{i\theta_k},\\
\rho_k=\cos\tfrac{\omega_{1,k}}2\cos\tfrac{\omega_{2,k}}2\in[0,1],\qquad
\theta_k=-\tfrac{\omega_{1,k}+\omega_{2,k}}2 .
\end{gathered}
\label{eq:half-pixel-kernel-multiplier}
\end{equation}
Thus $S$ is a mode-wise contraction, with $s_k=\rho_ka_k$.
$\theta_k$ is the relative phase of the two models.
\end{proposition}
\begin{proof}
Factor the symbol $\tfrac14\sum_{r,t\in\{0,1\}}e^{-i(r\omega_1+t\omega_2)}$
and use $1+e^{-i\omega}=2\cos(\omega/2)e^{-i\omega/2}$; on the signed
frequency interval each cosine lies in $[0,1]$.
\end{proof}

The common Fourier basis, $\mu_k=1$, and the proposition verify
\cref{ass:standing}. Since $\gamma_k=1$ only at the zero frequency, the
nondegeneracy condition \eqref{eq:active-model-nondegeneracy} holds as soon
as the data excite a nonzero frequency with $\nu_k\ne0$. Then
\cref{thm:main-spectral-general,thm:phase-aware-threshold} apply with
$s_k=\rho_ka_k$. The next result shows where the half-pixel phase controls
the threshold and where it cannot.

\begin{corollary}[Thresholds for the half-pixel pair]
\label{cor:half-pixel-thresholds}
Let $k$ be an active frequency with $\nu_k\ne0$ and $0<\rho_k<1$.
\begin{enumerate}[label=(\roman*)]
\item If $16\rho_k\bigl(1-\cos\frac{\omega_{1,k}+\omega_{2,k}}2\bigr)\ge(1-\rho_k)^2$,
then $c_k$ decreases strictly on $(0,\infty)$ and $\beta_k^\theta=0$.
\item If $\omega_{1,k}+\omega_{2,k}=0$, then $\theta_k=0$,
$\beta_k^\theta=\beta_{k,+}=a_k^2\bigl(\rho_k+(1+\rho_k)\sqrt{\rho_k}\bigr)$,
and $k$ satisfies the sign-change condition.
\end{enumerate}
Consequently,
\(
\lambda^\theta\ge
\max\bigl\{a_k^2\bigl(\rho_k+(1+\rho_k)\sqrt{\rho_k}\bigr):
k\in\mathcal K,\ \omega_{1,k}+\omega_{2,k}=0,\ \nu_k\ne0,\ 0<\rho_k<1\bigr\},
\)
and, by \cref{thm:phase-aware-threshold}(iv), no threshold that depends
only on the spectra and the active set can lie below this bound.
\end{corollary}
\begin{proof}
(i) is \eqref{eq:phase-dominance} divided by $a_k^2$. For (ii), $t_k=0$, so
\cref{thm:phase-aware-threshold}(ii) gives $\beta_k^\theta=\beta_{k,+}$, a
simple zero of $P_k=d_k(\beta-p_k)P_{0,k}$ at which $P_k$ changes sign.
\end{proof}

At low frequency, $1-\rho_k\approx|\omega_k|^2/8$ and
$1-\cos\theta_k\approx(\omega_{1,k}+\omega_{2,k})^2/8$, with
$|\omega_k|^2=\omega_{1,k}^2+\omega_{2,k}^2$. Condition (i) then holds
outside a thin wedge
$|\omega_{1,k}+\omega_{2,k}|\lesssim|\omega_k|^2/(8\sqrt2)$ around the
anti-diagonal. The half-pixel shift thus acts mainly as a phase shift,
except on frequencies whose wave vector is orthogonal to the shift
direction $(1,1)$, where it is a pure attenuation. For a normalized kernel
on a square grid, $a_k\to1$ and $\rho_k\to1$ at the lowest anti-diagonal
frequencies $\pm(1,-1)$, so the bound tends to $3$ as the grid is refined.
On a rectangular grid the exactly anti-diagonal frequencies are the
multiples of $\pm(n_1,-n_2)/\gcd(n_1,n_2)$, which may lie at high frequency.
The threshold can then be attained at nearly anti-diagonal modes, which
violate \eqref{eq:phase-dominance} without being phase-aligned
(\cref{app:rect-grid}).

\subsection{What the half-pixel discrepancy measures}
\label{subsec:roughness}

In general, $D(\lambda)$ measures how sensitive the Tikhonov reconstruction
is to the prescribed model perturbation. For half-pixel perturbations a
more tangible reading is available. Consider the Nyquist-truncated
pure-phase half-pixel shift $S_\varphi=F\diag(\gamma^\varphi)F^*$,
\begin{equation}
\gamma_k^\varphi=\begin{cases}
e^{-i(\omega_{1,k}+\omega_{2,k})/2},&\omega_{1,k},\omega_{2,k}\in(-\pi,\pi),\\
0,&\text{on the Nyquist lines}.
\end{cases}
\label{eq:exp-purephase}
\end{equation}
The truncation preserves conjugate symmetry and hence real outputs. This
operator is a mode-wise contraction with modal magnitudes $0$ or $1$, so its
data-independent threshold is zero. Off the Nyquist lines it has the same
phase $\theta_k$ as the bilinear multiplier.

\begin{proposition}[Half-pixel discrepancy as directional roughness]
\label{prop:roughness}
Let $L=I_n$ and $\widetilde A=F\diag(\nu)F^*$.
\begin{enumerate}[label=(\roman*)]
\item \emph{Pure-phase shift (exact).} For $S=S_\varphi \forall\lambda>0$, $\widetilde x_\lambda^S=S_\varphi^*\widetilde x_\lambda$, hence
$D(\lambda)=\|(I-S_\varphi^*)\widetilde x_\lambda\|_2$. Off the Nyquist
lines, the modulus of $I-S_\varphi^*$ is
$|1-\overline{\gamma_k^\varphi}|=2|\sin((\omega_{1,k}+\omega_{2,k})/4)|$,
which equals $|\omega_{1,k}+\omega_{2,k}|/2$ to leading order: a
half-pixel difference along the shift direction $(1,1)$.
\item \emph{Bilinear interpolation (leading order).} Let
$c_k^\varphi:=2(1-\cos\theta_k)\phi_k^2$ be the modal coefficient of the
pure-phase pair. For every active mode with $a_k>0$ and $\rho_k>0$ and every
$\lambda>0$,
\[
|c_k(\lambda)-c_k^\varphi(\lambda)|
\le\frac{1-\rho_k}{\rho_k^2}\Bigl(\frac{1-\rho_k}{\rho_k^2}
+2(1-\cos\theta_k)\Bigr)\phi_k(\lambda)^2 .
\]
\end{enumerate}
\end{proposition}
\begin{proof}[Proof of (i)]
If $|\gamma_k^\varphi|=1$, the second filter
$\overline{\nu_k\gamma_k^\varphi}/(|\nu_k|^2+\lambda)$ is
$\overline{\gamma_k^\varphi}$ times the first; if $\gamma_k^\varphi=0$, both
sides vanish. The symbol follows from $|1-e^{i\psi}|=2|\sin(\psi/2)|$.
(ii) is proved in \cref{app:roughness-proof}.
\end{proof}

At low frequency, the bound in (ii) is $O(|\omega_k|^4)\phi_k^2$, whereas $c_k^\varphi$ is of order $(\omega_{1,k}+\omega_{2,k})^2\phi_k^2$. Thus, along any fixed frequency direction not orthogonal to $(1,1)$, the bilinear and pure-phase coefficients agree to leading order as the frequency tends to zero.
This need not hold near the anti-diagonal or at high frequencies, whose contributions also depend on the data weights. For the pure-phase pair, the target prescribes exactly the directional
roughness $\|(I-S_\varphi^*)\widetilde x_\lambda\|_2/\|b\|_2$.
For the bilinear pair, this interpretation holds to leading order at
low frequency in the directions described above. It gives $\tau$ a
concrete meaning: an acceptable level of variation along the shift
direction relative to the data. This resembles norm-constrained
(Ivanov-type) regularization, where a prescribed bound is related to
a penalized formulation through a Lagrange multiplier
\cite{ivanov2013theory}. The formulations differ: SSDM matches a
roughness measure while retaining the penalty $\|x\|_2^2$, normalizes
the target by $\|b\|_2$, and transfers the selected
parameter to the original operator.

\section{Numerical experiments}
\label{sec:numerical-experiments}

This section presents three deblurring experiments to test four claims listed below. Experiment~1 is a clean test of (C1) where the theory applies exactly. Experiment~2 introduces the mismatch that motivates transfer to test (C1) and (C2). Experiment~3 repeats on natural photographs whose blur model is only approximately right to test (C1) to (C3). We separately analyse (C4).

\emph{(C1) Shape and certificate}
(\cref{sec:spectral-theory,sec:controller}).
Under the stated assumptions, the surrogate disagreement $D_b(\lambda)$
has a strict maximum at $\lambda=0$, and its square decays as
$\lambda^{-2}$. For the bilinear half-pixel pair, the data-independent
threshold guarantees decrease only beyond $\lambda^\theta$, which lies
between about $1.6$ and $3$ in the tested settings. The finite,
data-dependent certificate can establish strict decrease below this
threshold, making each target in the certified range correspond to a
unique weight. For the pure-phase pair, the threshold is zero, so the
spectral theory already guarantees strict decrease for every $\lambda>0$.

\emph{(C2) Transfer} (\cref{sec:framework}, Step~4). The weight chosen with
the cheap surrogate can be used with the original operator. However, theory does not guarantee a good reconstruction. 

\emph{(C3) The target sets the selected level}. Selection is certified given a target, which determines the regularization level. The reconstruction quality depends on the setting (forward model, noise, penalty, boundary treatment) and is not covered by the theory. 

\emph{(C4) Root finding is cheap} with an certified interval; certificate dominates cost. (\cref{sec:exp-root-solving}).
% \begin{enumerate}[label=(C\arabic*),leftmargin=*]
% \item \emph{Shape and certificate} (\cref{sec:spectral-theory,sec:controller}).
% For the pair of surrogate models, the disagreement $D_b(\lambda)$ between the
% two reconstructions is largest at $\lambda=0$, and its square decays like
% $\lambda^{-2}$. It is guaranteed to decrease only beyond a threshold
% $\lambda^\theta$ of order $1$ to $3$. In the practical range, decrease is a property of the data, which a finite certificate can check. Where the certificate holds, every target in range is met at exactly one weight.
% \item \emph{Transfer} (\cref{sec:framework}, Step~4). The weight chosen with
% the cheap surrogate is used with the original operator. The theory does not
% promise that this gives a good reconstruction. The experiments measure what
% transfer repairs and what it does not.
% \item \emph{The target sets the selected level.} Given a target, selection is
% certified. The target determines which regularization level is selected, but
% whether that level gives an accurate reconstruction depends on the whole
% setting (forward model, noise, penalty, boundary treatment) and is not covered
% by the theory. The experiments record, retrospectively, how a good target behaves.
% \item \emph{Root finding is cheap} once the certificate holds. The certificate
% is the dominant cost of selection (\cref{sec:exp-root-solving}).
% \end{enumerate}
\paragraph{Evaluation.}
The input to parameter selection consists of the data, the two surrogate models and a prescribed target $\tau$~(\cref{sec:framework}). Ground truth is used \emph{post hoc} to evaluate reconstruction accuracy through the \emph{relative error}
\(
E=\|x-x_{\mathrm{true}}\|_2/\|x_{\mathrm{true}}\|_2
\),
and the \emph{numerical oracle}, which estimates the smallest error achievable by varying the regularization weight for the original operator. Its weight and error are denoted by $\lambda_{\mathrm{or}}$ and $E_{\mathrm{or}}$, and the corresponding \emph{oracle-associated target} is $\tau_{\mathrm{or}}:=D_b(\lambda_{\mathrm{or}})$. Results are reported mainly as $E/E_{\mathrm{or}}$, with $1$ denoting the oracle reference.

In Experiments~1 and~2, the oracle is obtained from $61$ logarithmically spaced weights followed by scalar-search refinement. In Experiment~3, its error is estimated by fitting a parabola through the three stored errors around the grid minimum. This estimate is not itself an evaluated reconstruction (\cref{app:exp3-setup}). Experiments~1 and~2 use targets fixed in advance. Experiment~3 additionally explores smaller targets after inspecting $\tau_{\mathrm{or}}$. Ground truth therefore informs this exploratory sweep, although it is never an input to the parameter-selection algorithm. Throughout, \emph{numerically verified} means that the certificate's conditions pass in double-precision arithmetic; no interval-arithmetic validation is performed.

\subsection{Common setup}
\label{sec:exp-setup}

\begin{figure}[t]
    \centering
    \includegraphics[width=\textwidth]{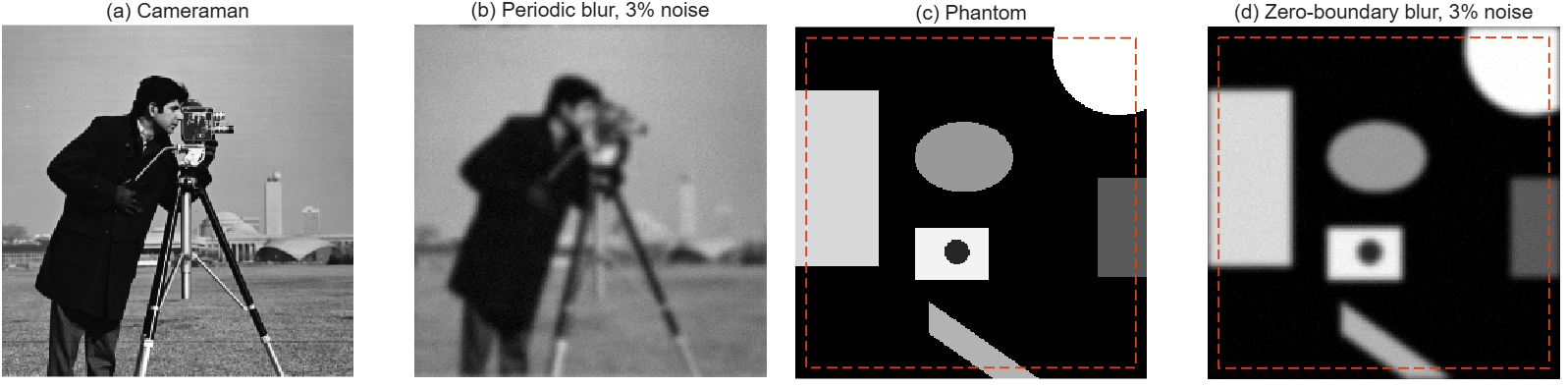}
    \caption{Ground truths and data. Experiment~1: (a)~Cameraman and
    (b)~its periodic blur with $3\%$ noise. Experiment~2: (c)~the
    piecewise-constant phantom, whose objects intersect all four boundaries,
    and (d)~its zero-boundary blur with $3\%$ noise. The dashed frame marks
    the outer eight-pixel strip over which boundary RMSE is computed.
    Display range $[0,1]$.}
    \label{fig:ground-truths}
\end{figure}

\paragraph{Operators and data.}
In Experiments~1 and~2, images have size $256\times256$ and are
column-vectorized ($n=256^2$). The PSF is the normalized separable Gaussian
$p=gg^\top$ with
\[
g_t=\frac{\exp(-t^2/(2\sigma^2))}{\sum_{s=-8}^{8}\exp(-s^2/(2\sigma^2))},
\qquad t=-8,\ldots,8,\qquad \sigma=2.5\text{ pixels},
\]
and the surrogate $\widetilde A=p*$ is periodic convolution
(\cref{sec:problem-setup}). Gaussian noise is scaled to $3\%$ of the
noise-free data norm. We write $x_B(\lambda;b)=(B^*B+\lambda I)^{-1}B^*b$ for the Tikhonov reconstruction with forward matrix $B$ and data $b$, $L=I_n$.
\Cref{fig:ground-truths} shows the ground truths and data.

\paragraph{Selection signal and targets.}
The primary perturbation is the bilinear half-pixel operator of
\cref{prop:half-pixel-multiplier}. For comparison we also use the
pure-phase shift \eqref{eq:exp-purephase}, for which the discrepancy is
exactly a directional roughness (\cref{prop:roughness}). Selection uses
\begin{equation}
D_b(\lambda)=\frac{\|x_{\widetilde A}(\lambda;b)
-x_{\widetilde AS}(\lambda;b)\|_2}{\|b\|_2}.
\label{eq:exp-Db}
\end{equation}
The targets are prescribed. The different target values are  used to test the method across a range of regularization levels. They are not recommendations (\cref{sec:target-choice}).
% : $\tau=0.05$ in Experiment~1 and
% $\tau\in\{0.05,0.075,0.10,0.125\}$ in Experiment~2. In Experiment~3 we start
% with the same four values and then add smaller ones, after seeing that the
% oracle-associated target is much lower there (\cref{sec:real-data}). 
% The targets

\paragraph{Certificate and root.}
The certificate of \cref{prop:interval-certificate} is applied to $\Phi_b$
on $[10^{-10},10^3]$ with $6000$ logarithmically spaced cells. Where it
holds, the root $\lambda^\star$ of $D_b=\tau$ is unique on this interval. We
compute it with the bracketing solver \texttt{fzero} in $\log\lambda$
(\texttt{TolX}$=10^{-10}$), started from the cell of a $241$-point grid
that contains the sampled crossing. Starting from the whole certified
interval gives the same roots to within a relative difference of
$1.3\times10^{-5}$. The cell count, about $460$ cells per decade, is a
numerical setting and not a tuned quantity
(\cref{prop:interval-certificate} guarantees success for a fine enough
partition where the derivative is strictly negative).

% \paragraph{Accuracy.}
% Accuracy is $E$ and $E/E_{\mathrm{or}}$ as defined above; $1.05$ means $5\%$
% above the best error found on the grid. These are evaluation quantities, not
% selection rules.

\subsection{Experiment 1: does the theory hold where it should?}
\label{sec:exp-spectral}

\begin{table}[t]
\centering
\caption{Experiment 1: principal results. The tail slope is that of
$\Phi_b$ versus $\lambda$ for $\lambda\ge10^2$ (predicted: $-2$).
Reconstruction errors are at a weight within the stopping band
$|D_b-\tau|\le2.5\times10^{-4}$ of the root.}
\label{tab:exp1-findings}
\small
\begin{tabular}{@{}lrr@{}}
\toprule
Quantity & Bilinear & Pure-phase\\
\midrule
Relative difference, direct versus spectral evaluation & $9.0\times10^{-15}$ & $1.6\times10^{-14}$\\
Fitted tail slope of $\Phi_b$ & $-1.9950$ & $-1.9950$\\
Closed-form threshold $\lambda^\ddagger$ & $2.98859$ & $0$\\
Phase-aware threshold $\lambda^\theta$ & $2.97722$ & $0$\\
Root $\lambda^\star$ at $\tau=0.05$ & $0.0055749$ & $0.0051891$\\
$\Phi_b'<0$ verified on $[10^{-10},10^3]$ & yes & yes\\
Reconstruction error $E$ ($E/E_{\mathrm{or}}$) & $0.13121$ ($1.0014$) & $0.13136$ ($1.0026$)\\
\bottomrule
\end{tabular}
\end{table}

We blur the Cameraman image with the periodic surrogate itself, add noise, and
run SSDM with the same surrogate. This is purposefully an inverse crime as the model that made the data is known. We aim to test the formulae and the certificate in the setting where they apply exactly. 

\paragraph{Design.}
The ground truth is Cameraman scaled to $[0,1]$, and
$b=\widetilde Ax_{\mathrm{true}}+\eta$. Selection and reconstruction thus use
the forward model that generated the data. The target is $\tau=0.05$. The
claim under test is (C1): that the direct and the spectral evaluation of
$D_b$ agree, that the tail decays as predicted, and that the certificate
passes although the root lies far below the data-independent threshold.

\paragraph{Findings.}
The direct and spectral evaluations of the discrepancy agree to machine precision, and the fitted tail slope matches the predicted $\lambda^{-2}$ decay of the squared discrepancy (\cref{tab:exp1-findings}). For the bilinear perturbation, the selected weight is approximately $\lambda^\theta/534$, where $\lambda^\theta=2.97722$. This threshold is attained at the anti-diagonal modes $\pm(1,-1)$ and is optimal among data-independent thresholds (\cref{cor:half-pixel-thresholds,app:exp-spectral-diagnostics}). Establishing monotonicity at the selected weight therefore requires the data-dependent analysis.

For both perturbations, the certificate verifies strict decrease throughout $[10^{-10},10^3]$, establishing uniqueness of the target root on this interval. The largest computed cell bound is approximately $-7.2\times10^{-13}$. Its small magnitude reflects the small derivative in the far tail. The bound closely tracks the largest sampled derivative, $-7.33\times10^{-13}$; moreover, at the sampled points, increasing modal contributions offset less than $2\%$ of the decreasing contributions (rise ratio below $0.02$). These diagnostics indicate that cancellation between modal contributions is mild.

At the bilinear root, the phase term in \eqref{eq:phase-decomposition} accounts for $98.7\%$ of $\Phi_b$, so the discrepancy is dominated by the relative phase of the two models. For both perturbations, the reconstruction error is within $0.3\%$ of the numerical-oracle error, $E_{\mathrm{or}}=0.13103$, compared with $E=0.16038$ for the observation. Thus, the prescribed target happens to be near-optimal in this matched setting. This is consistent with (C1). Next we examine whether reconstruction accuracy persists under model mismatch.

\subsection{Experiment 2: what does transfer fix?}
\label{sec:exp-transfer}

Now the data are blurred with zero values assumed outside the image (a
zero-boundary model), which the periodic surrogate cannot reproduce near the
edges. The phantom has objects touching all four edges, so that the edges
matter. We select $\lambda$ with the surrogate pair, hand it to the original
blur, and ask what is gained. As a control we repeat the experiment with data
that the surrogate itself generated, which removes the mismatch.

\paragraph{Design.}
A piecewise-constant phantom intersects all four image boundaries
(\cref{fig:ground-truths}(c,d)). The original operator $A$ uses same-size
convolution with $p$ and zero exterior values, while the surrogate
$\widetilde A$ uses periodic boundaries. We compare zero-boundary data
$b_0=Ax_{\mathrm{true}}+\eta_0$ with a matched control
$b_{\mathrm P}=\widetilde Ax_{\mathrm{true}}+\eta_{\mathrm P}$, using
the same noise array scaled to $3\%$ of each noise-free data norm.
The signal-dependent mismatch is $7.48\%$, more than twice the noise level.

Selection uses the bilinear surrogate pair in both cases
(M-B and Z-B). For zero-boundary data, the selected weight is transferred
to $A$; the surrogate reconstruction at the same weight is retained to
isolate the effect of transfer. Reconstructions are evaluated at the
target roots to separate target effects from stopping tolerances.
Original-model solves use preconditioned conjugate gradients, checked
against a direct solve. Boundary and interior RMSE are measured on the
outer eight-pixel strip and its complement, respectively.

\begin{table}[t]
    \centering
    \caption{Experiment~2: root-based sensitivity for the bilinear pair.
    $E/E_{\mathrm{or}}$ uses the numerical oracle of each data model (matched:
    $\lambda=0.0050119$, $E=0.148096$; zero-boundary: $\lambda=0.0063096$,
    $E=0.134630$). Boundary RMSE is given for the transferred and the
    surrogate reconstruction at the same $\lambda$.}
    \label{tab:exp2-sensitivity}
    \footnotesize\setlength{\tabcolsep}{4pt}
    \begin{tabular}{@{}crrr|rrrrr@{}}
        \toprule
        & \multicolumn{3}{c|}{Matched periodic (M-B)} & \multicolumn{5}{c}{Zero-boundary transfer (Z-B)} \\
        $\tau$ & $\lambda^\star$ & $E$ & $E/E_{\mathrm{or}}$ & $\lambda^\star$ & $E$ & $E/E_{\mathrm{or}}$
               & \multicolumn{2}{c}{Bdry.\ RMSE: transf.\,/\,surr.} \\
        \midrule
        0.050 & 0.0563463 & 0.17350 & 1.172 & 0.115177 & 0.19630 & 1.458 & 0.1021 & 0.1471 \\
        0.075 & 0.0023340 & 0.15159 & 1.024 & 0.0231257 & 0.14327 & 1.064 & 0.0619 & 0.1593 \\
        0.100 & 0.00065250 & 0.18082 & 1.221 & 0.0078334 & 0.13506 & 1.003 & 0.0538 & 0.2049 \\
        0.125 & 0.00033943 & 0.21607 & 1.459 & 0.0037837 & 0.13542 & 1.006 & 0.0528 & 0.2541 \\
        \bottomrule
    \end{tabular}
\end{table}

\begin{figure}[tb]
    \centering
    \includegraphics[width=\textwidth]{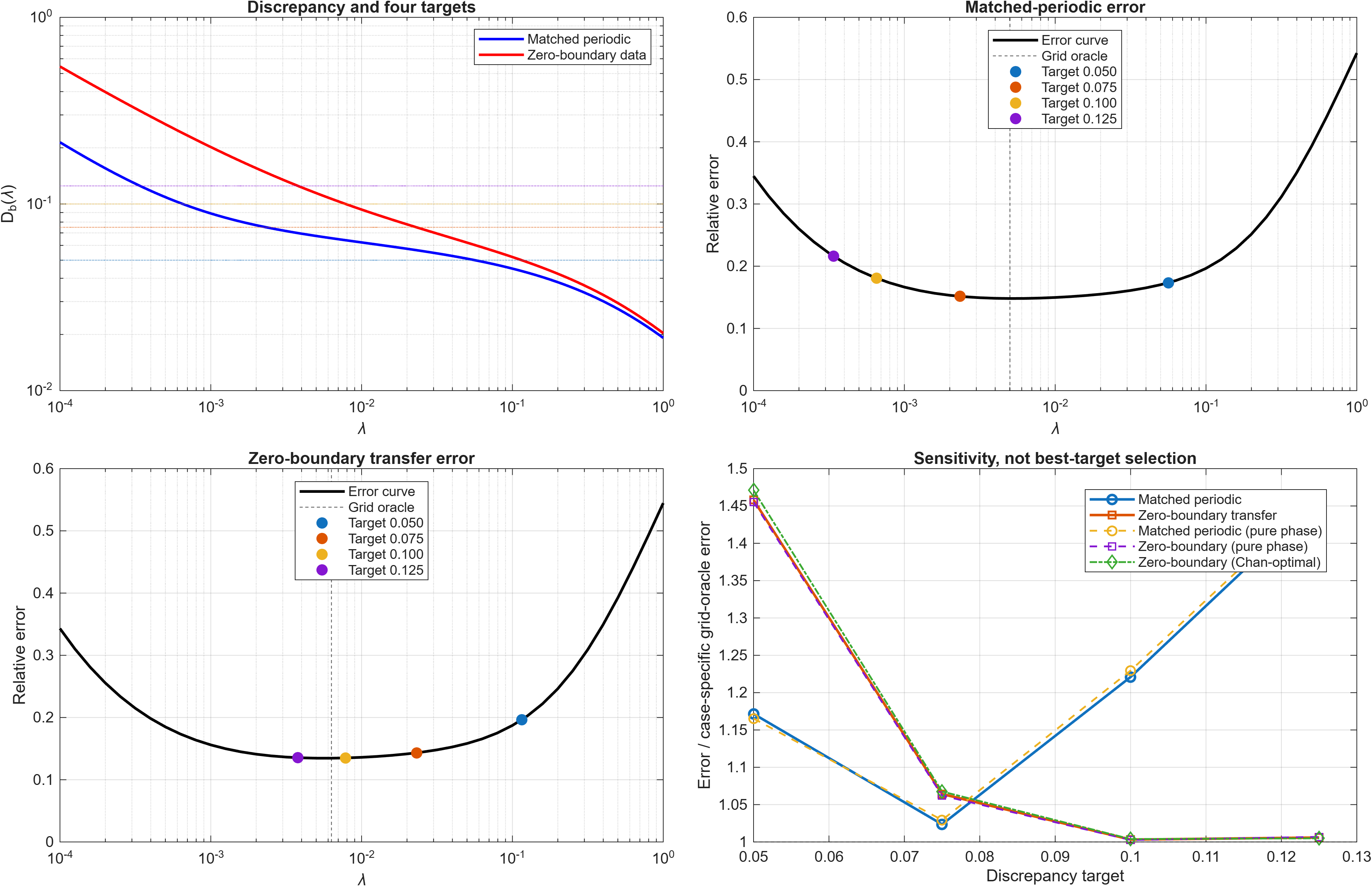}
    \caption{Experiment~2: (a)~$D_b$ for matched and zero-boundary data with the four targets;
    (b,c)~relative-error curves (evaluation only) with the roots and numerical oracles;
    (d)~oracle-normalized error across targets for all five cases.}
    \label{fig:exp2-paired}
\end{figure}

\paragraph{The certificate holds under mismatch.}
We certify five cases: the bilinear and pure-phase pairs on matched and
zero-boundary data, and the bilinear pair with the Frobenius-nearest
circulant surrogate~\cite{chan1988optimal} on zero-boundary data. For every
case the certificate is verified on all of $[10^{-10},10^3]$, so each $D_b$
curve decreases strictly there and each target root is unique. For the
bilinear pair this is not implied by the spectral theory: all roots lie
below $0.12$, whereas the phase-aware thresholds of the periodic and
nearest-circulant surrogates are $2.97722$ and $2.88710$. The zero-boundary
curve lies above the matched one, so its roots are larger
(\cref{fig:exp2-paired}(a)); this ordering is observed, not predicted, since
changing the data changes both the Fourier energies and the normalization.

\paragraph{Transfer improves the boundary; the target sets the level.}
At the four selected weights, transfer to $A$ reduces boundary RMSE
by $31$--$79\%$ relative to the periodic reconstruction
(\cref{tab:exp2-sensitivity}), but does not correct an unsuitable
regularization level. For zero-boundary data, $\tau=0.10$ and $0.125$
give errors within $0.3\%$ and $0.6\%$ of the numerical oracle.
The same targets yield smaller weights and amplified noise in the
matched control, with error ratios of $1.22$ and $1.46$.
Among the tested targets, $\tau=0.075$ gives the best compromise,
within $2.4\%$ and $6.4\%$ of the respective oracles.
These retrospective comparisons separate the benefits of boundary
correction from the setting-dependent choice of target.

\paragraph{Surrogate and original-operator roots differ.}
Separability of $A$ allows us to evaluate
\(
D_b^A(\lambda)=\|x_A(\lambda;b_0)-x_{AS_0}(\lambda;b_0)\|_2/\|b_0\|_2,
\)
where $S_0$ is the zero-exterior bilinear shift
(\cref{app:exp2-solvers}). This diagnostic is not part of SSDM. It
measures how target matching changes when the surrogate pair is replaced
by $(A,AS_0)$. Each target has one sampled crossing, without a uniqueness
certificate (\cref{tab:exp2-root-shift}). The roots differ by about $9\%$
at $\tau=0.05$, while at the larger targets SSDM selects weights
$2.5$--$8$ times larger. Nevertheless, its transferred reconstructions
can be more accurate: at $\tau=0.10$ and $0.125$, they are near the
oracle, whereas the original-operator roots give error ratios of
$1.13$ and $1.40$. The ordering reverses at $\tau=0.075$.
Thus reconstruction accuracy does not require the surrogate root to
approximate the original-operator root. This example also illustrates
that a target's meaning depends on the model pair
(\cref{sec:target-choice}); it concerns one dataset and one choice of
$S_0$.

\begin{table}[tb]
    \centering
    \caption{Experiment~2, case Z-B: the root of $D_b=\tau$ for the surrogate
    pair (SSDM) and of $D_b^A=\tau$ for the original pair $(A,AS_0)$, and the
    transferred error $E/E_{\mathrm{or}}$ at each (numerical oracle as in
    \cref{tab:exp2-sensitivity}). $D_b^A(\lambda_{\mathrm{SSDM}})$ is the
    original-operator discrepancy at the surrogate root.}
    \label{tab:exp2-root-shift}
    \small
    \begin{tabular}{@{}crrrrrr@{}}
        \toprule
        $\tau$ & $\lambda_{\mathrm{SSDM}}$ & $\lambda_A$ & $\lambda_A/\lambda_{\mathrm{SSDM}}$ & $D_b^A(\lambda_{\mathrm{SSDM}})$ & $E/E_{\mathrm{or}}$ at $\lambda_{\mathrm{SSDM}}$ & $E/E_{\mathrm{or}}$ at $\lambda_A$ \\
        \midrule
        0.050 & 0.115177 & 0.125616 & 1.09 & 0.0510 & 1.458 & 1.504 \\
        0.075 & 0.0231257 & 0.00914445 & 0.40 & 0.0670 & 1.064 & 1.007 \\
        0.100 & 0.0078334 & 0.00115730 & 0.15 & 0.0763 & 1.003 & 1.128 \\
        0.125 & 0.0037837 & 0.00046587 & 0.12 & 0.0832 & 1.006 & 1.403 \\
        \bottomrule
    \end{tabular}
\end{table}

\paragraph{Reference rules under model mismatch.}
We compare SSDM with generalized cross-validation (GCV) and the
discrepancy principle (DP), using the same surrogate and transferring
the selected weights to the original operator. DP uses the true noise
norm to solve
\(
\|\widetilde A x_{\widetilde A}(\lambda;b)-b\|_2
=\zeta\|\eta\|_2.
\)
For zero-boundary data, we also evaluate DP with the original-operator
residual. The reference rules are competitive under matched conditions
(\cref{tab:exp2-reference}), but boundary mismatch causes surrogate GCV
and DP with $\zeta=1$ and $1.1$ to underregularize severely, with errors
about $292$, $111$ and $57$ times the numerical-oracle error,
respectively. The surrogate residual contains model error
($\delta=0.07477$) exceeding the noise level ($\varepsilon=0.03$),
so matching it to noise alone encourages fitting that mismatch.

Inflating the DP level restores accuracy: $\zeta=1.5$ gives
$E/E_{\mathrm{or}}=1.002$, while the quadrature estimate
$\zeta=\sqrt{1+(\delta/\varepsilon)^2}=2.69$ gives $1.213$.
Similar behavior occurs with the Frobenius-nearest circulant surrogate
(\cref{tab:exp2-reference}). Both DP and SSDM admit near-oracle choices
in retrospective sweeps, but require different external inputs:
a residual level accounting for noise and mismatch for DP, and a
discrepancy target for SSDM. Here, $\delta$ uses ground truth and
$\zeta=1.5$ was chosen retrospectively. These results illustrate
sensitivity to model mismatch on one image and noise realization;
they do not establish a ranking of the methods.

\begin{table}[tb]
    \centering
    \caption{Experiment~2: reference rules on the same surrogate, with
    transfer to $A$ for zero-boundary data. SSDM entries give the range over
    the four roots of \cref{tab:exp2-sensitivity}. DP: discrepancy principle
    with the true noise norm. For $\zeta=2.69=\sqrt{1+(\delta/\varepsilon)^2}$ (a
    quadrature estimate) the level also uses the mismatch $\delta=0.07477$,
    which is computed from $x_{\mathrm{true}}$, so this row is an
    oracle-informed benchmark. The $\zeta=1.5$ value does not use $\delta$ but
    was chosen with the sweep results in view. Neither row was computed for the
    matched data.}
    \label{tab:exp2-reference}
    \footnotesize\setlength{\tabcolsep}{4pt}
    \begin{tabular}{@{}llrrrr@{}}
        \toprule
        & & \multicolumn{2}{c}{Matched (M)} & \multicolumn{2}{c}{Zero-boundary (Z)} \\
        \cmidrule(lr){3-4}\cmidrule(lr){5-6}
        Rule & Uses & $\lambda$ & $E/E_{\mathrm{or}}$ & $\lambda$ & $E/E_{\mathrm{or}}$ \\
        \midrule
        SSDM, $\tau\in\{0.05,\ldots,0.125\}$ & $\tau$ & $3.4\times10^{-4}$--$0.056$ & $1.024$--$1.459$ & $0.0038$--$0.115$ & $1.003$--$1.458$ \\
        Surrogate GCV & --- & $0.001757$ & $1.046$ & $2.559\times10^{-8}$ & $292$ \\
        DP, $\zeta=1$, surrogate residual & $\|\eta\|_2$ & $0.006574$ & $1.001$ & $1.482\times10^{-7}$ & $111$ \\
        DP, $\zeta=1.1$, surrogate residual & $\|\eta\|_2$ & $0.01346$ & $1.023$ & $4.201\times10^{-7}$ & $56.6$ \\
        DP, $\zeta=1.5$, surrogate residual & $\|\eta\|_2$, $\zeta$ & --- & --- & $0.004338$ & $1.002$ \\
        DP, $\zeta=2.69$, surrogate residual & $\|\eta\|_2$, $\delta$ & --- & --- & $0.05895$ & $1.213$ \\
        DP, $\zeta=1$, original residual & $\|\eta\|_2$, $A$ & --- & --- & $0.006531$ & $1.000$ \\
        DP, $\zeta=1.1$, original residual & $\|\eta\|_2$, $A$ & --- & --- & $0.01326$ & $1.023$ \\
        Numerical oracle & $x_{\mathrm{true}}$ & $0.0050119$ & $1$ & $0.0063096$ & $1$ \\
        \bottomrule
    \end{tabular}
\end{table}

% \paragraph{Findings.}
% For (C1), the certificate holds in all five cases although every root lies far
% below the data-independent threshold. For (C2), transfer repairs what the
% surrogate gets wrong at the boundary and leaves the regularization level where the
% target put it. For (C3), the same target is good in one data model and poor in
% the other. The surrogate pair, selecting on its own, still reaches an accurate
% reconstruction for several targets, but not because its root approximates the
% root of the original operator.

\subsection{Experiment 3: do natural photographs behave the same way?}
\label{sec:real-data}

Now we ask whether the behavior established on the smaller test images persists for natural photographs when both the blur model and the surrogate’s boundary treatment are imperfect. We also examine what the monotonicity guarantees mean for reconstruction accuracy. The blur is generated computationally, so the sharp photographs remain available for evaluation.

\begin{figure}[tb]
    \centering
    \includegraphics[width=0.92\textwidth]{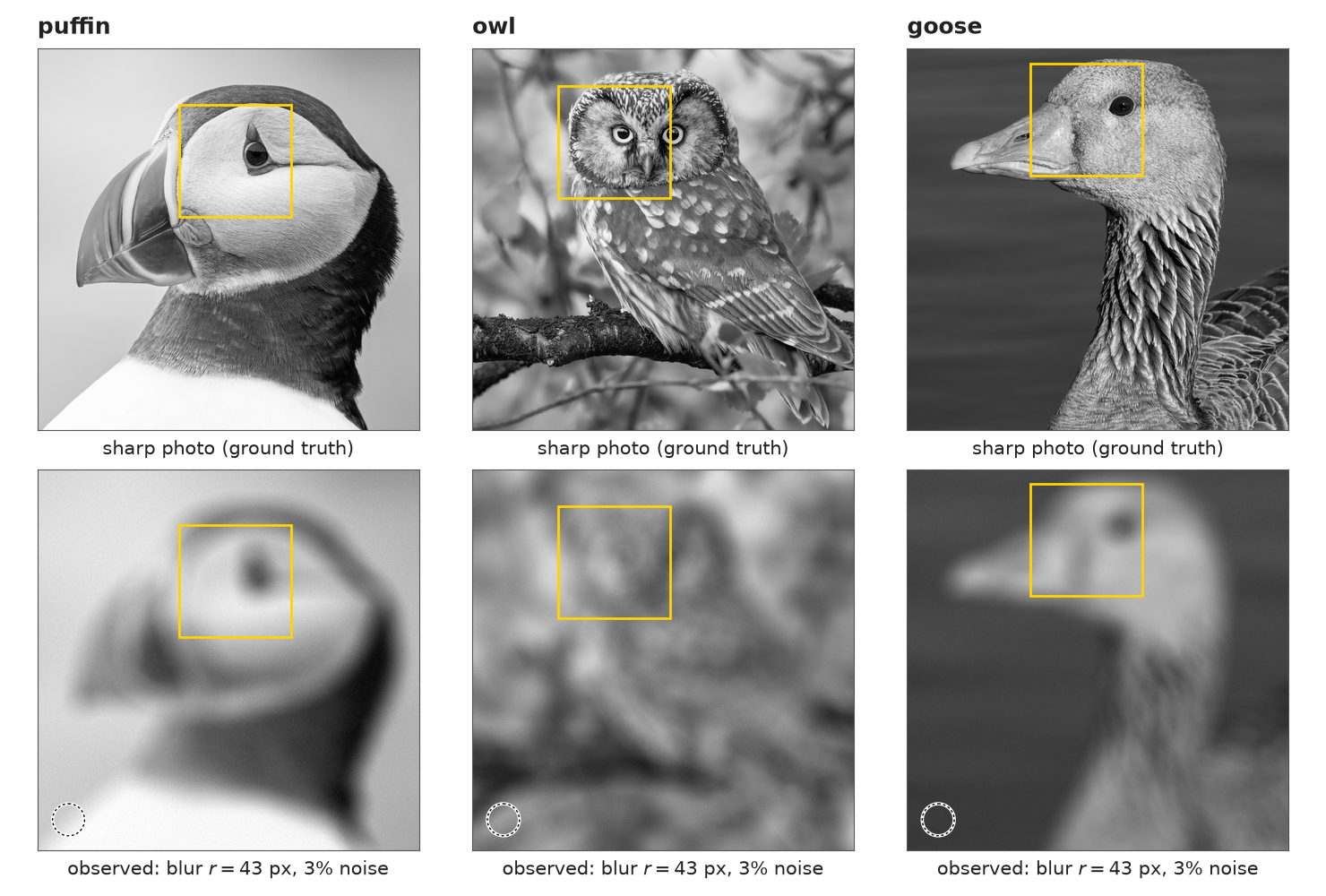}
    \caption{Experiment~3 data (the $3\%$ noise variant). Top: the sharp
    photographs, which serve as ground truth for scoring only. Bottom: the
    observations, blurred by a disc of radius $43$~px and corrupted by noise;
    the circle in the lower left shows the size of the disc. The yellow frame
    marks the region shown in the zoom of \cref{fig:bird-recon}.}
    \label{fig:bird-data}
\end{figure}

\paragraph{Data.}
Three grayscale photographs---a puffin, an owl and a goose
(\cref{fig:bird-data})---are prepared as $1110\times1110$ fields.
The central $1024\times1024$ window is the ground truth
$x_{\mathrm{true}}$. A normalized disc PSF of radius $43$ pixels
blurs the full field, and only its valid $1024\times1024$ convolution
is recorded. Thus the data include contributions from outside the
reconstruction window without imposing an artificial boundary condition.
Gaussian noise is added at $3\%$ and $1\%$ of the noise-free data norm,
using the same noise pattern rescaled. The sharp photographs are used
for evaluation and the retrospective target study. All figures use
display range $[0,1]$.

\paragraph{Models and reconstruction.}
\label{sec:real-data-setup}
The main setting assumes a radius-$39$ disc PSF, four pixels smaller
than the true PSF. A radius-$43$ control isolates boundary effects.
For either radius, the surrogate $\widetilde A$ implements reflexive
blur by periodic convolution on the $2048\times2048$ mirror tile of
the recorded data. The second model $\widetilde AS$ applies the
periodic half-pixel shift on this tile. Both are BCCB and share
the Fourier basis required by \cref{sec:spectral-theory,sec:controller}, though the shifted model is not a reflexive-boundary model and does not preserve mirror symmetry. Transfer uses the truncated convolution operator $A$, with unknown on the recorded window extended by the assumed PSF radius, and the valid convolution on that window as output. Its normal
equations are solved by preconditioned conjugate gradients, as in
Experiment~2. With the exact-radius control, $A$ matches the
data-generating operator, while the surrogate retains its boundary
mismatch.

% \begin{itemize}
%   \item The surrogate $\widetilde A$ is the disc blur with reflexive boundary
%   conditions, realized as a BCCB matrix on the $2\times2$ mirror tile of the
%   data ($2048\times2048$) and diagonalized by the FFT. The second model applies
%   the periodic half-pixel shift on this tile. It does not preserve the mirror
%   symmetry, so $\widetilde AS$ is a BCCB model on the tile and not a
%   reflexive-boundary model. The theory of
%   \cref{sec:spectral-theory,sec:controller} requires only the common Fourier
%   basis on the tile.
%   \item The original operator $A$ is the truncated, non-circulant model: the
%   unknown lives on the window extended by the PSF radius, the data only on the
%   window, and $A$ is the \texttt{'valid'} part of the convolution with the PSF. Its normal
%   equations are solved by preconditioned conjugate gradients, as in
%   Experiment~2.
%   \item The model disc has radius $39$ (\emph{main} setting), so the PSF is
%   wrong by four pixels. In the \emph{control} setting it has radius $43$, the
%   true one. Then $A$ reproduces the data up to the noise, and the only defect
%   left is the boundary model of the surrogate.
% \end{itemize}

\paragraph{Magnitude of the model mismatch.}
For these photographs, the radius-$39$ PSF introduces only a mild data
mismatch: after accounting for noise in quadrature, its contribution is
approximately $1.2$--$1.7\%$ of the data norm. It is therefore smaller than
the noise at $3\%$ but larger at $1\%$, motivating the second noise level.
The discrepancy between the original operator and the reflexive surrogate
is also small, about $0.9$--$1.9\%$, compared with $7.5\%$ in Experiment~2.
These tests thus examine mild PSF and boundary mismatch. Reconstruction
errors are measured against the original photographs; allowing an affine
gain and offset does not materially change the conclusions
(\cref{app:exp3-setup}).

\paragraph{Settings, targets and comparisons.}
The experiment has $3$ photographs $\times$ $2$ noise levels $\times$ $2$ model
radii, that is, twelve configurations, each run with both perturbations. The
targets start with the four values of Experiment~2. The first runs showed that
these lie far above $\tau_{\mathrm{or}}\approx0.01$, and we then added smaller
targets ($0.0075$, $0.01$, $0.015$, $0.02$, $0.03$, and $0.04$ for the
control), \emph{after} seeing the oracle-associated targets. The sweep over
small targets is therefore exploratory. Selection, certificate (for the
radius-$39$ configurations), and root finding are as in
\cref{sec:exp-setup}. We compare with the observation itself, with the
surrogate weights selected by GCV and by the L-curve criterion (maximum
curvature), both transferred to $A$ like SSDM, and with the numerical oracle. All metrics are computed on the $1024\times1024$ window.

\paragraph{The certificate holds on photographs.}

For the six radius-$39$ configurations and both perturbations, the certificate
is verified on the whole interval $[10^{-10},10^3]$, with largest cell bound
between $-9.5\times10^{-15}$ and $-8.0\times10^{-15}$ and rise ratio $0.030$
(bilinear pair). The phase-aware threshold is $\lambda^\theta=2.979$, close to
the value $3$ of the closed form, because the tile is square
(\cref{cor:half-pixel-thresholds}), while all selected roots are at most
$0.01$. Decrease in the selection range is a property of the data, and
here it holds. In the control configurations the certificate was not computed.
There we checked only that $D_b$ decreases along the $241$-point grid on
$[10^{-8},10^3]$ and that its analytic derivative is negative on a $1201$-point
grid, which is evidence and not a certificate.

\begin{figure}[tb]
    \centering
    \includegraphics[width=\textwidth]{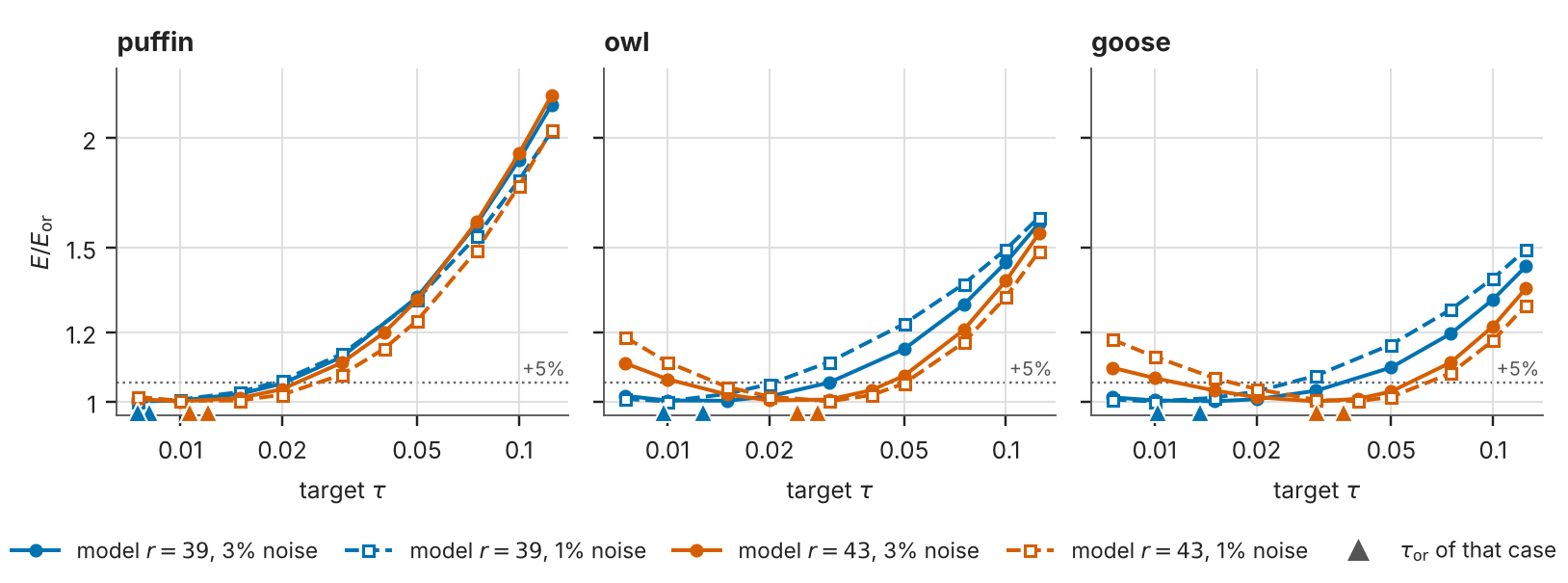}
    \caption{Experiment~3: accuracy against target for the bilinear pair.
    $E/E_{\mathrm{or}}=1$ is the error of the numerical oracle for that
    configuration, and the dotted line marks $5\%$ above it. Colors distinguish the model radius, line styles the noise level. Triangles mark the
    oracle-associated target $\tau_{\mathrm{or}}$ of each curve. The four
    targets of Experiment~2 ($0.05$ to $0.125$) lie on the right of each panel.
    The pure-phase pair gives curves within $0.012$ of these at $\tau=0.015$
    and $0.02$.}
    \label{fig:bird-sweep}
\end{figure}

\paragraph{Targets from the Gaussian experiments are mostly too large.}
At the targets inherited from Experiment~2 ($\tau\ge0.05$), the transferred
reconstructions have errors between $1.01$ and $2.24$ times the oracle
error (\cref{fig:bird-sweep}). Few are within $5\%$ of the oracle, and
some are worse than the blurred observations. The oracle-associated
targets instead range from $0.0075$ to $0.036$ for the bilinear pair,
with similar values for the pure-phase pair. The larger targets select
insufficient regularization and amplify noise. This illustrates the
distinction between target matching and target suitability: the roots
are found, with uniqueness guaranteed where certified, but the inherited
targets do not transfer reliably to these photographs.

\paragraph{Useful target ranges depend on the setting.}
Every configuration admits a range of targets giving reconstruction
errors within $5\%$ of the oracle error (\cref{fig:bird-sweep,tab:bird-results}). These ranges shift with the model and noise level: the oracle-associated targets are consistently lower for the radius-$39$ model than for the exact model, and respond differently to changes in noise. Despite this variation, common targets perform well across the tested configurations. For the bilinear pair, $\tau=0.015$ and $0.02$ give errors within $6.4\%$ and $5.7\%$ of the oracle, respectively, in all twelve configurations, with similar results for the pure-phase pair. These targets were chosen retrospectively: they demonstrate a useful shared range in the tested settings, rather than a target-selection rule available in advance.

\begin{figure}[H]
    \centering
    \includegraphics[width=\textwidth]{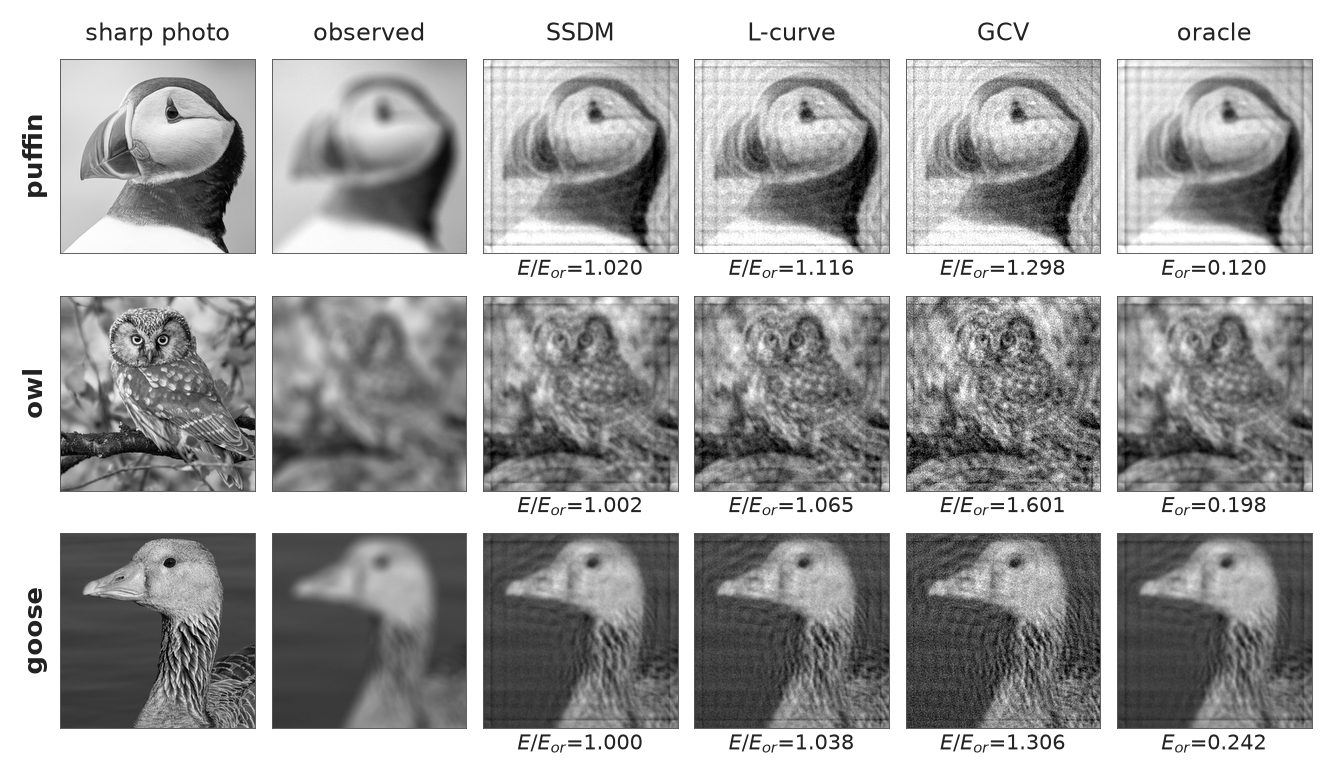}\\\vspace{5mm}
    \includegraphics[width=\textwidth]{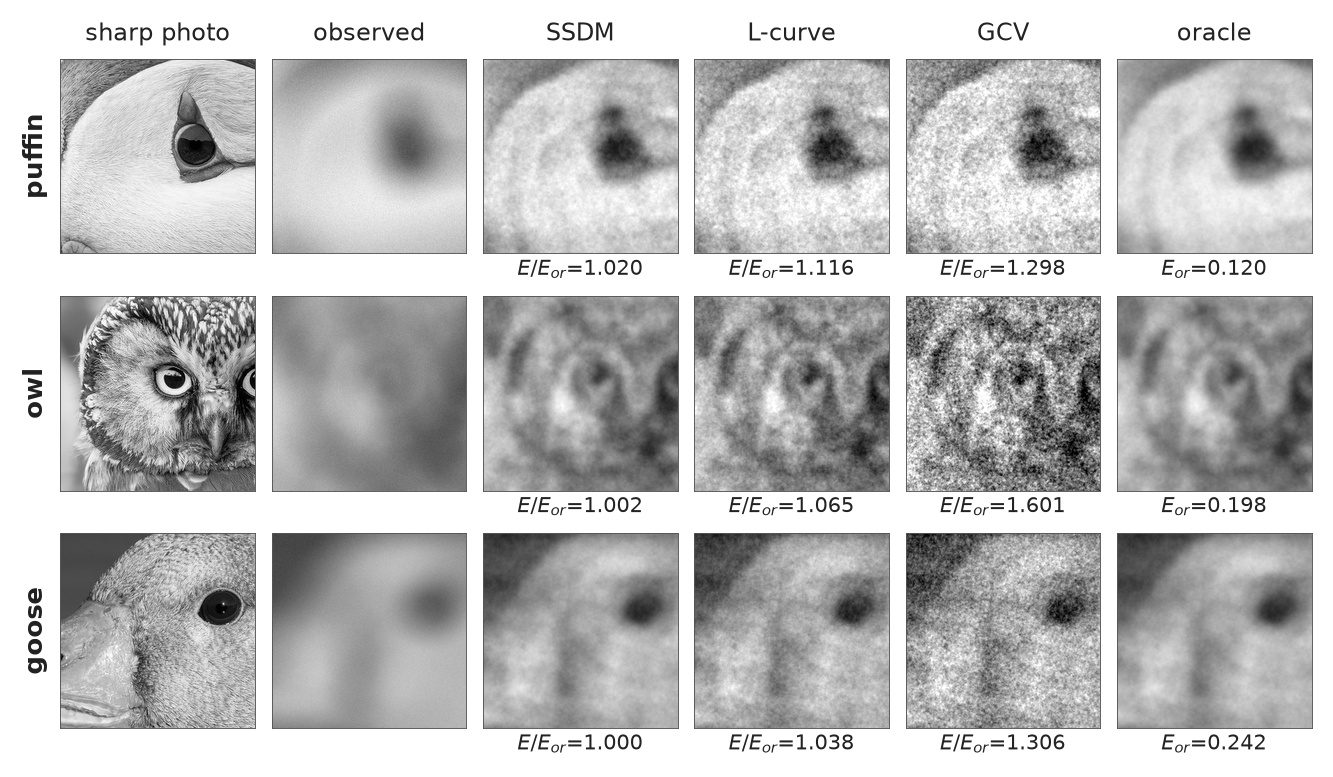}\vspace{5mm}
    \caption{Experiment~3: reconstructions for the radius-$39$ model and $3\%$ noise.
    Columns: sharp photograph, observation, SSDM (bilinear pair, $\tau=0.015$,
    transferred to $A$), transferred L-curve and GCV choices, and the numerical
    oracle, with $E/E_{\mathrm{or}}$ below each reconstruction (below the oracle, its error $E_{\mathrm{or}}$). The oracle is a benchmark and not a method: the reconstruction with the original operator at the weight that has the smallest error against the original photograph, which is known here only because the data are synthetic. Top: whole images at half resolution. Bottom: the $300\times300$ region marked in \cref{fig:bird-data}. SSDM is close to the oracle in $E/E_{\mathrm{or}}$, and both retain visible artifacts. The L-curve and GCV choices carry visibly more noise.}
    \label{fig:bird-recon}
\end{figure}

\paragraph{Reference rules.}
At $\tau=0.015$, SSDM gives lower reconstruction error than the
transferred L-curve choice in nine of twelve configurations and than
GCV in all twelve (\cref{tab:bird-results}). The L-curve and GCV errors
range from $1.00$--$1.38$ and $1.17$--$2.06$ times the oracle error,
respectively. The observation itself is only $1.13$--$1.54$ times the
oracle error, illustrating that global relative error captures only part of the benefit of deblurring; differences in edges and texture are visible in \cref{fig:bird-recon}. Since $\tau=0.015$ was chosen retrospectively using the oracle, whereas the reference rules use no ground truth, these comparisons do not establish a ranking. Instead, they show that a suitably chosen target can yield competitive reconstruction accuracy, while leaving its prospective calibration unresolved.

\begin{table}[H]
\centering
\caption{Experiment~3: oracle and transferred SSDM reconstructions for the
twelve configurations (bilinear pair), with errors against the original photo. $E_{\mathrm{or}}$: numerical-oracle error;
$\tau_{\mathrm{or}}$: oracle-associated target; window: range of targets with
$E/E_{\mathrm{or}}\le1.05$, from linear interpolation of the sweep in
$\log\tau$ (``$\le$'' means that $E/E_{\mathrm{or}}\le1.05$ holds down to the
smallest target tested, $0.0075$). The remaining columns give
$E/E_{\mathrm{or}}$ for SSDM at $\tau=0.015$ and $0.02$, the transferred
L-curve and GCV choices, and the blurred observation. $r$: model disc radius
(true radius $43$).}
\label{tab:bird-results}
\footnotesize\setlength{\tabcolsep}{3.4pt}
\begin{tabular}{@{}llrrrcrrrrr@{}}
\toprule
& noise & $r$ & $E_{\mathrm{or}}$ & $\tau_{\mathrm{or}}$ & window
& \multicolumn{2}{c}{SSDM $\tau=$} & L-curve & GCV & obs. \\
\cmidrule(lr){7-8}
& & & & & & $0.015$ & $0.02$ & & & \\
\midrule
\multirow{4}{*}{Puffin} & 3\% & 39 & 0.120 & 0.0081 & $\le$0.020 & 1.020 & 1.048 & 1.116 & 1.298 & 1.139 \\
 & 1\% & 39 & 0.117 & 0.0075 & $\le$0.019 & 1.025 & 1.057 & 1.382 & 1.469 & 1.145 \\
 & 3\% & 43 & 0.114 & 0.0107 & $\le$0.022 & 1.009 & 1.033 & 1.093 & 1.228 & 1.201 \\
 & 1\% & 43 & 0.107 & 0.0121 & $\le$0.025 & 1.003 & 1.018 & 1.252 & 1.169 & 1.248 \\
\midrule
\multirow{4}{*}{Owl} & 3\% & 39 & 0.198 & 0.0127 & $\le$0.030 & 1.002 & 1.013 & 1.065 & 1.601 & 1.222 \\
 & 1\% & 39 & 0.196 & 0.0097 & $\le$0.020 & 1.020 & 1.047 & 1.384 & 2.058 & 1.229 \\
 & 3\% & 43 & 0.173 & 0.0242 & 0.011--0.045 & 1.018 & 1.003 & 1.011 & 1.294 & 1.402 \\
 & 1\% & 43 & 0.156 & 0.0278 & 0.014--0.050 & 1.039 & 1.013 & 1.114 & 1.454 & 1.544 \\
\midrule
\multirow{4}{*}{Goose} & 3\% & 39 & 0.242 & 0.0137 & $\le$0.035 & 1.000 & 1.006 & 1.038 & 1.306 & 1.128 \\
 & 1\% & 39 & 0.239 & 0.0102 & $\le$0.025 & 1.010 & 1.027 & 1.265 & 1.849 & 1.133 \\
 & 3\% & 43 & 0.215 & 0.0302 & 0.012--0.056 & 1.029 & 1.011 & 1.001 & 1.167 & 1.270 \\
 & 1\% & 43 & 0.193 & 0.0362 & 0.017--0.063 & 1.064 & 1.033 & 1.048 & 1.213 & 1.402 \\
\bottomrule
\end{tabular}
\end{table}

\begin{table}[b]
\centering
\caption{Variations of the oracle-associated target $\tau_{\mathrm{or}}$ and
the oracle parameter $\lambda_{\mathrm{or}}^{\mathrm{ref}}$ in the paired
noise-level study of \cref{app:tau-noise} (medians over $20$ realizations,
noise levels $0.3$--$10\%$; \cref{tab:tau-noise-medians}). Ratios are
max/min over the listed values. These values are not used to select targets.}
\label{tab:target-vs-parameter}
\small
\begin{tabular}{@{}lrr@{}}
\toprule
Variation & $\tau_{\mathrm{or}}$ & $\lambda_{\mathrm{or}}^{\mathrm{ref}}$ \\
\midrule
Across noise levels, matched data & $0.0612$--$0.0798$ ($1.30\times$) & $1.15\times10^{-4}$--$0.0253$ ($220\times$) \\
Across noise levels, zero-boundary data & $0.0772$--$0.437$ ($5.7\times$) & $1.39\times10^{-4}$--$0.0254$ ($183\times$) \\
Across boundary models, fixed noise & $1.26$--$5.5\times$ & $1.00$--$1.21\times$ \\
\bottomrule
\end{tabular}
\end{table}

\paragraph{The weight need not match the oracle weight for the accuracy to match.}
At $\tau=0.015$ the selected weight is between $0.40$ and $3.6$ times the oracle
weight, and the error is still within $6.4\%$. The error curve $E(\lambda)$ is
flat near its minimum for these photographs, so a factor of two or three in
$\lambda$ costs a few percent. This is why the useful ranges of targets are wide, and also why a comparison of weights alone would overstate the differences.

\paragraph{Transfer has image-dependent effects.}
At $\tau=0.015$, transfer to the original operator changes the
reconstruction error by $-5.5\%$ to $+1.5\%$ for the owl and goose,
consistent with the small mismatch between the original operator and
the surrogate. For the puffin, however, transfer increases the error
by $8$--$17\%$, including an increase in interior RMSE. Before transfer,
the surrogate reconstruction even outperforms the original-operator
oracle in all four settings, illustrating that this oracle is optimal
only within its own reconstruction family. The cause of the puffin's
deterioration remains unresolved. These results do not contradict
(C1) or (C3), which concern the surrogate pair; they show that improved
accuracy after transfer is a separate property, not guaranteed by the
theory.

\subsection{Root finding and computational cost}
\label{sec:exp-root-solving}

On a certified interval that brackets the target, the root is unique.
Bracketing in $\log\lambda$ reached the stopping tolerance without
failure in all $22$ root problems of Experiments~1 and~2.
Experiment~3 used the same root-finding procedure, although its
radius-$43$ controls were checked only by sampling rather than interval
certification. These results support the procedure in
\cref{alg:ssdm}: certify an interval, then bracket the root.

After an $O(n\log n)$ data transform, interval certification costs
$O(Jn)$ for $J$ cells and $n$ active modes, while each subsequent
evaluation of $D_b$ is an $O(n)$ spectral sum. Certification dominates
the root search and serves all targets in the certified range for that
image. Its guarantee does not transfer to new data. 

The final reconstruction requires a solve with the original operator
($44$--$223$ PCG iterations in Experiment~3). No wall-clock times were
recorded, and no computational advantage over other parameter rules
is established. In particular, the separable operator of Experiment~2
provides a controlled transfer test rather than evidence of an advantage
when the original spectrum is inaccessible
(\cref{app:exp2-solvers}).

\subsection{Target dependence and interpretation}
\label{sec:target-choice}

The oracle-associated target $\tau_{\mathrm{or}}$ retrospectively
identifies the discrepancy level associated with the most accurate
reconstruction. In the paired noise-level study of Experiment~2,
it varies substantially less than the oracle parameter when the model
is fixed, especially for matched data. Across boundary models, however,
the targets can differ considerably even when the oracle parameters
nearly coincide
(\cref{app:tau-noise,tab:target-vs-parameter}).

Experiment~3 shows similar stability under changes in noise and greater
sensitivity to changes in the PSF. Its oracle-associated targets lie
below those used in the Gaussian experiments, although common targets
$\tau=0.015$ and $0.02$ give errors within $7\%$ of the oracle in all
twelve photographic configurations (\cref{tab:bird-results}).
These retrospective observations establish neither general
transferability nor that choosing $\tau$ is easier than choosing
$\lambda$.

\Cref{prop:roughness} gives the target an interpretation as directional
variation relative to the data: exactly for the pure-phase pair and
to leading order at low frequency away from the orthogonal direction
for the bilinear pair. Its suitability nevertheless depends on the
image, noise, model mismatch and perturbation. Small disagreement alone
does not establish accuracy, since strong regularization can remove
useful features from both reconstructions.

\paragraph{Relation to the theory.}
\label{sec:exp-summary}
The certificate inequalities held in floating-point
arithmetic. For the bilinear pair, this supports strict decrease and
unique target matching in practical ranges well below the
data-independent threshold. For the pure-phase pair, the spectral
theory already guarantees decrease for every positive parameter.
These results concern the surrogate pair. Transfer reduced boundary
errors in Experiment~2, but changed little for two photographs and
increased the puffin's error by $8$--$17\%$. The experiments support the  selection mechanism while separating it from the setting-dependent questions of target calibration and reconstruction accuracy. The \emph{post-hoc} comparisons do not establish superiority over other parameter rules.

\section{Conclusion}
\label{sec:conclusion}

Double-model parameter selection approaches regularization through the
agreement of two reconstructions of the same data. In the SSDM framework,
a structured surrogate and a perturbation sharing its spectral basis make
this comparison analytically tractable. Given the model pair, penalty and
comparison measure, a prescribed disagreement target is the sole scalar input to parameter selection, whose evaluation requires neither ground truth nor a known noise level.

The central contribution is to establish when this target determines a
unique regularization parameter. The spectral analysis shows that
individual modal contributions can be nonmonotone, so uniqueness cannot
be assumed. The phase-aware threshold guarantees eventual decrease
independently of the data, while the finite certificate verifies
decrease on intervals relevant to selection. For half-pixel deblurring,
the connection to directional roughness further explains what the
disagreement measures. The experiments verified monotonicity in practical
selection ranges well below the universal threshold, while showing that
transfer to the original operator does not invariably improve accuracy.

A prescribed target separates a modelling decision from its numerical
realization. It expresses how much disagreement between the two
representations is acceptable, and the regularization parameter is then
determined by matching that requirement. This offers an interpretable
way to specify regularization through sensitivity to a controlled model
change, without comparing the reconstruction to an unavailable true
image. Choosing an appropriate level of sensitivity remains a separate
question: uniqueness makes the requirement unambiguous, but does not
ensure that it preserves the features needed for a particular task.

The results therefore provide a foundation for both target calibration
and accuracy analysis. Future work should investigate how targets can
be prescribed from application requirements or inferred from data, and
when they transfer across images, noise levels and model pairs. Relating
disagreement to reconstruction error or task-relevant features would
clarify when agreement is informative and when it merely reflects
excessive smoothing. Further directions include designing perturbations
that probe meaningful model uncertainty and extending the monotonicity
analysis beyond simultaneously diagonalizable quadratic problems.

% ------------------------------------------------------------
% Bibliography
% ------------------------------------------------------------
\bibliographystyle{plainnat}
\bibliography{references}

% ------------------------------------------------------------
% Appendices
% ------------------------------------------------------------
\clearpage
\appendix
\pdfbookmark[1]{Appendices}{appendices}
\begin{center}
{\Large\bfseries Appendices\par}
\end{center}
\vspace{0.5em}

\section{Supporting analysis}
\label{app:controller-details}

This appendix supplies the remaining proofs, a multiple-root example,
and an extension beyond contractive perturbations.

\subsection{Example: a target with three roots}
\label{app:controller-multiple-crossings}

A target may be attained more than once. Consider the scalar pair
$\widetilde A=1$, $S=1/2$, $L=1$ with data $b=1$. It satisfies
the contraction, positive-penalty, and nondegeneracy hypotheses, and
\[
D(\lambda)=\left|\frac{1}{1+\lambda}-\frac{1/2}{1/4+\lambda}\right|,\qquad
D(0)=1,\quad D(1/2)=0,\quad D(1)=\frac1{10},\quad
\lim_{\lambda\to\infty}D(\lambda)=0 .
\]
The intermediate value theorem gives a crossing of $\tau=1/20$ in each of
$(0,1/2)$, $(1/2,1)$, and $(1,\infty)$. By
\cref{lem:scalar-filter-difference}, $D$ decreases on $(0,1/2)$, increases
on $(1/2,\lambda^\ddagger)$, and decreases on $(\lambda^\ddagger,\infty)$,
with $\lambda^\ddagger=\tfrac12+\tfrac{3}{2\sqrt2}>1$, which equals
$\lambda^\theta$ because the relative phase is zero. Since $D(1)>\tau$,
there is no crossing in $[1,\lambda^\ddagger]$. Thus there are exactly
three positive roots. Only the last lies in the
data-independent tail; uniqueness on a practical interval requires a
certificate such as \cref{prop:interval-certificate}.

\subsection{Proof of \texorpdfstring{\cref{thm:phase-aware-threshold}}{Theorem 3.8}(ii)--(iv)}
\label{app:phase-aware-proofs}

\paragraph{Part (ii): phase dependence.}
If $s_k=0<a_k$, then $P_k(\beta)=a_k^2\beta^3>0$ for $\beta>0$.
If $a_k=s_k>0$, then $P_k(\beta)=2t_ka_k^2(a_k^2+\beta)N_k(\beta)$, which is
positive for $\beta>0$ if $t_k>0$ and vanishes identically if $t_k=0$. If
$a_k=0$, then $P_k\equiv0$. In all these cases $\beta_k^\theta=0$. Now let
$a_k>s_k>0$ and write
\[
P_k=P_k^{(0)}+t_kE_k,\qquad
P_k^{(0)}:=d_k(\beta-p_k)P_{0,k},\qquad
E_k:=p_k(a_k^2+s_k^2+2\beta)N_k>0.
\]
By \eqref{eq:beta-pm},
$P_k^{(0)}>0$ on $(0,p_k)\cup(\beta_{k,+},\infty)$ and $P_k^{(0)}<0$ on
$(p_k,\beta_{k,+})$. If $t_k=0$, then $\beta_k^\theta=\beta_{k,+}$. If
$t_k>0$, then $P_k>0$ on $[\beta_{k,+},\infty)$, and continuity gives
$\beta_k^\theta<\beta_{k,+}$. For $0\le t<t'$, the inequality
$P_k^{(0)}+t'E_k\ge P_k^{(0)}+tE_k$ shows that the set in
\eqref{eq:beta-theta-def} shrinks as $t_k$ grows, which gives
monotonicity, and $\lambda^\theta\le\lambda^\ddagger$ follows.

\paragraph{Part (iii): phase dominance and sharpness.}
Suppress $k$, let $a>s>0$, write $\sigma:=a^2+s^2>2p$, and set
\[
G(\beta):=(\sigma+2\beta)N(\beta)+16(\beta-p)P_0(\beta).
\]
Regarding $G$ as
a function of $\beta$, $p$, and $\sigma$, with $N=p^2+\sigma\beta+\beta^2$
and $P_0=\beta^2-2p\beta-p(\sigma+p)$, we have for $\sigma\ge2p$
\[
\frac{\partial G}{\partial\sigma}=p^2+2\sigma\beta+3\beta^2-16p(\beta-p)
\ge3(\beta-2p)^2+5p^2>0,
\qquad
G\big|_{\sigma=2p}=2(\beta+p)(3\beta-5p)^2\ge0 .
\]
Hence $G>0$ for $\beta>0$ when $\sigma>2p$. Condition
\eqref{eq:phase-dominance} reads $tp\ge d/16$, so
$P(\beta)\ge d(\beta-p)P_0+\frac d{16}(\sigma+2\beta)N=\frac d{16}G(\beta)>0$
and $\beta^\theta=0$.

To prove sharpness, let $c>16$, $a=1$, $s=\rho\in(0,1)$,
and $t=(1-\rho)^2/(c\rho)$, so that $c\,as\,t=(a-s)^2$. At
$\beta=5\rho/3$,
\[
\frac{P(\beta)}{d}=\frac{G(\beta)}{16}
-\Bigl(\frac1{16}-\frac1c\Bigr)(\sigma+2\beta)N(\beta)
\;\longrightarrow\;
-\Bigl(\frac1{16}-\frac1c\Bigr)\frac{16}3\cdot\frac{64}9<0
\quad(\rho\uparrow1),
\]
since $G=0$ at $\sigma=2p$, $\beta=5p/3$. Thus $\beta^\theta\ge5\rho/3>0$
for $\rho$ close to $1$, although $c\,as\,t\ge(a-s)^2$ (and $t\le2$, so $t$
is realized by some $\theta$).

\paragraph{Part (iv): sharpness of the threshold.}
By \cref{lem:exact-modal-derivative},
$c_k'(\lambda_1)=-2\mu_kP_k(\beta_1)/N_k(\beta_1)^3>0$, so
\[
\Phi_\omega'(\lambda_1)=c_k'(\lambda_1)\omega_k+\sum_{j\ne k}c_j'(\lambda_1)\omega_j>0
\quad\text{whenever}\quad
\omega_k>\frac{\sum_{j\ne k}|c_j'(\lambda_1)|\omega_j}{c_k'(\lambda_1)}.
\]
If $\bar\lambda<\lambda^\theta=\beta_k^\theta/\mu_k$ for a maximizing mode
satisfying \eqref{eq:modal-sign-change}, there is
$\beta_1\in(\bar\lambda\mu_k,\beta_k^\theta)$ with $P_k(\beta_1)<0$, and the
first part gives a weight vector for which $\Phi_\omega$ increases at
$\lambda_1=\beta_1/\mu_k>\bar\lambda$.

The two suprema in \eqref{eq:beta-theta-def} and \eqref{eq:modal-sign-change}
differ only when the largest positive zero of $P_k$ is one at which $P_k$
does not change sign. An isolated zero of the derivative does not preclude
strict monotonicity, so without the sign-change condition $\lambda^\theta$
need not be the smallest monotonicity threshold.

\subsection{Beyond contraction}
\label{rem:beyond-contraction}

Contraction makes the second model nonamplifying, but is not essential
to the spectral conclusions. For any simultaneously unitarily
diagonalizable pair with positive active penalties, $c_k$, $h_k^2$, the
phase term, $\beta_{k,+}$ and $P_k$ are symmetric in the two magnitudes.
Thus \cref{thm:main-spectral-general,thm:phase-aware-threshold} extend by
ordering the magnitudes mode by mode and replacing $\mathcal J$ with
$\{k\in\mathcal K:a_ks_k>0,\ a_k\ne s_k\}$. Modes with exactly one zero
magnitude decrease strictly for every $\lambda>0$.

\subsection{Proof of \cref{prop:interval-certificate}}
\label{app:certificate-proof}

\begin{proof}
\emph{Cell bounds.} Fix $j$ and $k$, let $\lambda\in[\bar\lambda_{j-1},\bar\lambda_j]$ and
$\beta=\mu_k\lambda\in[\beta_{k,j}^-,\beta_{k,j}^+]$. By
\cref{lem:exact-modal-derivative},
$c_k'(\lambda)=(2\mu_k/N_k(\beta)^3)(-P_k(\beta))$, where the factor
$2\mu_k/N_k(\beta)^3$ is positive and decreasing in $\beta$, and
$-P_k(\beta)\le-m_{k,j}$. If $m_{k,j}\ge0$, the nonpositive number
$-m_{k,j}$ is multiplied by a factor at least $2\mu_k/N_k(\beta_{k,j}^+)^3$;
if $m_{k,j}<0$, the positive number $-m_{k,j}$ is multiplied by a factor at
most $2\mu_k/N_k(\beta_{k,j}^-)^3$. In both cases $c_k'(\lambda)\le u_{k,j}$.
Multiplying by $|\widehat b_k|^2$ and summing gives $\Phi'\le U_j$.

\emph{Monotonicity and target matching.} If every $U_j<0$, then $\Phi'<0$ on $I$, so $\Phi\ge\Phi(\bar\lambda_J)$ on
$I$, and $\Phi(\bar\lambda_J)>0$: otherwise its strictly negative
derivative at that endpoint would contradict $\Phi\ge0$. Hence $D=\sqrt\Phi$ is positive and smooth on $I$,
with $D'=\Phi'/(2D)<0$, and continuity and strict decrease give the unique
root. If $\bar\lambda_J>\lambda^\theta$,
\cref{thm:phase-aware-threshold}(i) extends the decrease to
$[\bar\lambda_0,\infty)$. The normalized statements follow by dividing $\Phi$ and its
derivative bounds by $\|b\|_2^2$.

\emph{Converse under refinement.} Write $u_{k,j}=F_k(\bar\lambda_{j-1},\bar\lambda_j)$, where
$F_k(\underline\lambda,\overline\lambda)$ is the right side of
\eqref{eq:certificate-upper} for the interval
$[\underline\lambda,\overline\lambda]$. The minimum of $P_k$ over
$[\mu_k\underline\lambda,\mu_k\overline\lambda]$ depends continuously on the
endpoints, and both cases of \eqref{eq:certificate-upper} give $0$ when the
minimum is $0$, so $F_k$ is continuous on
$\{\bar\lambda_0\le\underline\lambda\le\overline\lambda\le\bar\lambda_J\}$,
with $F_k(\lambda,\lambda)=c_k'(\lambda)$. Uniform continuity gives
$\max_j\sup_{[\bar\lambda_{j-1},\bar\lambda_j]}|U_j-\Phi'|\to0$ as the mesh
tends to zero, and $\max_I\Phi'<0$ gives the claim.
\end{proof}

\paragraph{Computing the cubic.}
\label{app:certificate-computation}
In monomial form $P_k(\beta)=\sum_{i=0}^3\pi_{k,i}\beta^i$ with
\begin{equation}
\begin{aligned}
\pi_{k,3}&=d_k+2t_kp_k=|\nu_k(1-\gamma_k)|^2,\\
\pi_{k,2}&=3p_k\bigl(t_k(a_k^2+s_k^2)-d_k\bigr),\\
\pi_{k,1}&=p_k\bigl(t_k[(a_k^2+s_k^2)^2+2p_k^2]-d_k(a_k^2-p_k+s_k^2)\bigr),\\
\pi_{k,0}&=p_k^2\bigl(d_k(a_k^2+p_k+s_k^2)+t_kp_k(a_k^2+s_k^2)\bigr).
\end{aligned}
\label{eq:Pk-coefficients}
\end{equation}

The cell minimum is found by evaluating $P_k$ at the endpoints and at
the real zeros of $P_k'$ inside the cell. Approximate critical points
need not give a lower bound on the true minimum, and floating-point
errors can have either sign. Validated certification therefore requires
enclosing the critical points, or bounding the polynomial over the whole
cell, while propagating input and arithmetic uncertainty, for example
with outward-rounded interval arithmetic.

\subsection{Proof of \texorpdfstring{\cref{prop:roughness}}{Proposition 5.3}(ii)}
\label{app:roughness-proof}

Suppress $k$ and let $f(r):=ra/(r^2a^2+\lambda)$, so that $\phi=f(1)$
and $\phi^S=f(\rho)$. For $r\in[\rho,1]$,
\[
|f'(r)|
=\frac{a|\lambda-r^2a^2|}{(r^2a^2+\lambda)^2}
\le\frac{a}{r^2a^2+\lambda}
\le\frac{a}{\rho^2(a^2+\lambda)}
=\frac{\phi}{\rho^2}.
\]
Hence $|h|=|f(1)-f(\rho)|\le(1-\rho)\phi/\rho^2$. By
\eqref{eq:phase-decomposition}, $c-c^\varphi=h^2-2(1-\cos\theta)\phi h$.
The bound on $|h|$ proves the claim.

\section{FFT conventions for the periodic surrogate}
\label{app:fourier-implementation}

Let $\mathcal F$ be the unnormalized two-dimensional forward FFT and
$\mathcal F^{-1}$ its inverse, including the factor $1/n$, where
$n=n_1n_2$. Arrays are identified with their column-vectorizations. The
convention of \cref{sec:problem-setup} then gives
\[
F^*b=\frac{\mathcal F(b)}{\sqrt n},\qquad
Fy=\sqrt n\,\mathcal F^{-1}(y),\qquad
\nu=\mathcal F(p),\qquad \gamma=\mathcal F(d_v).
\]
Embed kernels by their periodic offsets: offset $(0,0)$ occupies the
first array entry, and negative offsets wrap around. Centered kernels
must first be rearranged into this layout. If interpolation offsets
coincide on a degenerate grid dimension, add their weights.

For $L=I_n$ and $\lambda>0$, reconstructions can be computed directly
from unnormalized FFTs as
\[
\widetilde x_\lambda
=\mathcal F^{-1}\!\left(
\frac{\overline\nu\,\mathcal F(b)}{|\nu|^2+\lambda}\right),\qquad
\widetilde x_\lambda^S
=\mathcal F^{-1}\!\left(
\frac{\overline{\nu\gamma}\,\mathcal F(b)}{|\nu\gamma|^2+\lambda}\right).
\]
Products and divisions inside the inverse FFT are entrywise. The
normalization factors cancel in these reconstructions, but Parseval
requires the factor $1/n$ in the discrepancy sums:
\[
\Phi(\lambda)=\frac1n\sum_k c_k(\lambda)|\mathcal F(b)_k|^2,
\qquad
\Phi_b(\lambda)=\sum_k c_k(\lambda)w_k,
\]
\[
w_k:=\frac{|\widehat b_k|^2}{\|b\|_2^2}
=\frac{|\mathcal F(b)_k|^2}{n\|b\|_2^2},\qquad
\sum_k w_k=1\quad(b\ne0).
\]
Omitting the factor $1/n$ from a discrepancy sum using unnormalized data
FFTs changes the signal and its target scale.

Use the same frequency ordering for $\nu$, $\gamma$ and the data
coefficients. The signed indices in \cref{prop:half-pixel-multiplier}
represent the Nyquist frequency by $-\pi$ when $n_j$ is even. The bilinear
multiplier vanishes exactly on either Nyquist line, where its phase is
irrelevant. Assigning these zeros by index avoids the small residual
from floating-point evaluation of $\cos(-\pi/2)$. This differs from
thresholding small nonzero blur eigenvalues, which changes the model
and can alter endpoint limits.

\section{Experiment 2: phantom, operators, and surrogates}
\label{app:exp2-solvers}
\paragraph{Phantom, operators, and data.}

The ground truth is a piecewise-constant $256\times256$ phantom with
objects crossing all four image boundaries
(\cref{fig:ground-truths}(c)). Objects are cropped at the field of view,
with zero exterior values. The original operator is same-size
zero-exterior convolution with $p$. By separability,
\[
\operatorname{mat}(A\operatorname{vec}X)=TXT^\top,\qquad A=T\otimes T,
\]
where $T\in\R^{256\times256}$ is symmetric Toeplitz and
$\operatorname{mat}$ inverts column-vectorization. Thus $A$ is BTTB and
is not diagonalized by $F$. The primary and matched-control data are
\[
b_0=Ax_{\mathrm{true}}+\eta_0,\qquad
b_{\mathrm P}=\widetilde Ax_{\mathrm{true}}+\eta_{\mathrm P}.
\]
Both noise vectors use the same Gaussian array, scaled separately to
$3\%$ relative noise. The signal-dependent mismatch is
$\norm{(A-\widetilde A)x_{\mathrm{true}}}_2/\norm{Ax_{\mathrm{true}}}_2
=0.07477$.

The phantom is sampled at the pixel centers $(X_j,Y_i)=((j-\tfrac12)/256,(i-\tfrac12)/256)$,
where $i$ indexes rows and $j$ columns. Starting from zero, the following sets are assigned in
order, each overwriting earlier values on overlaps:
$\{X<0.24,\ 0.18<Y<0.68\}\mapsto0.85$;
$\{(X-0.92)^2+(Y-0.06)^2<0.19^2\}\mapsto1$;
$\{|Y-0.56-0.70X|<0.045,\ X>0.38\}\mapsto0.70$;
$\{X>0.86,\ 0.43<Y<0.71\}\mapsto0.35$;
$\{((X-0.48)/0.14)^2+((Y-0.37)/0.10)^2<1\}\mapsto0.60$;
$\{0.34<X<0.55,\ 0.57<Y<0.72\}\mapsto0.95$;
$\{(X-0.46)^2+(Y-0.64)^2<0.035^2\}\mapsto0.15$.
The generating script is provided as accompanying material.

\paragraph{Transfer and references.}

For zero-boundary data, the selected weight is transferred to $A$ as in
\eqref{eq:framework-transfer}. PCG solves
$(A^*A+\lambda I_n)x=A^*b_0$ to relative tolerance $10^{-8}$, with at
most $500$ iterations and applied inverse preconditioner
$(\widetilde A^*\widetilde A+\lambda I_n)^{-1}$. The surrogate
reconstruction at the same weight is retained for comparison. For
matched data, both selection and reconstruction use $\widetilde A$.
No positivity constraint, clipping or edge tapering is applied.

For a direct reference solve, write
$T=U\operatorname{diag}(d)U^\top$, with $U\in\R^{256\times256}$
orthogonal; this $U$ is distinct from the common surrogate basis $Q=F$.
For the data array $B_{\mathrm{img}}:=\operatorname{mat}(b)$, set
$M:=dd^\top$ and $B_U:=U^\top B_{\mathrm{img}}U$. Then
\[
x_A(\lambda;b)=\operatorname{vec}\!\left[
U\left(\frac{M\odot B_U}{M\odot M+\lambda\mathbf1\mathbf1^\top}\right)U^\top
\right],
\]
where multiplication $\odot$ and the displayed division are entrywise, and
$\mathbf1\in\R^{256}$ is the vector of ones.
This direct formula gives the same unique Tikhonov solution. It is used
for numerical oracles, the transferred surrogate-GCV baseline, and solver
checks. The adjoint test gives $3.7\times10^{-17}$; the separable and
Toeplitz forms reproduce $A$ to $6.0\times10^{-16}$ and
$2.8\times10^{-16}$. Every PCG solve meets the residual tolerance and
agrees with the direct solution to within $2.6\times10^{-7}$.

Alongside $E$, PSNR and SSIM, we report boundary RMSE on the outer
eight-pixel strip and interior RMSE on its complement. The accessible
eigendecomposition makes this a controlled test of selection and
transfer, not evidence of a computational advantage when the original
operator's spectral decomposition is inaccessible. No timing claim is made.

\paragraph{Original-operator discrepancy (diagnostic).}
\label{app:exp2-original-discrepancy}
To compare selection signals from the surrogate and original operators
(\cref{tab:exp2-root-shift}), let $S_0$ be the zero-exterior version
of the bilinear shift, $(S_0X)_{ij}=\tfrac14(X_{ij}+X_{i-1,j}+X_{i,j-1}+X_{i-1,j-1})$
with $X=0$ outside the array. Then $S_0=T_{\mathrm s}\otimes T_{\mathrm s}$,
where $T_{\mathrm s}\in\R^{256\times256}$ is lower bidiagonal with entries
$\tfrac12$ on the diagonal and the first subdiagonal, and
$AS_0=K\otimes K$ with $K=TT_{\mathrm s}$. With the singular value
decomposition $K=U_K\operatorname{diag}(\varsigma)V_K^\top$, $M_K:=\varsigma\varsigma^\top$ and
$B_K:=U_K^\top B_{\mathrm{img}}U_K$,
\[
x_{AS_0}(\lambda;b)=\operatorname{vec}\!\left[
V_K\left(\frac{M_K\odot B_K}{M_K\odot M_K+\lambda\mathbf1\mathbf1^\top}\right)V_K^\top
\right],
\]
and $D_b^A(\lambda)=\|x_A(\lambda;b)-x_{AS_0}(\lambda;b)\|_2/\|b\|_2$, with
$x_A$ from the formula above. The roots $\lambda_A$ of $D_b^A=\tau$ are
located like the surrogate roots (sign changes on the $241$-point grid on
$[10^{-8},10^3]$, then \texttt{fzero} in $\log\lambda$).
Both operators are separable, so this diagnostic needs no iterative
solves. The identities $S_0=T_{\mathrm s}\otimes T_{\mathrm s}$ and
$A=T\otimes T$ hold to $10^{-16}$, and the surrogate roots reproduce
\cref{tab:exp2-sensitivity}. The choice of $S_0$ is not unique, and the
certificate of \cref{sec:controller} does not apply to $D_b^A$. Each
target has one sampled crossing; this is not a uniqueness guarantee.

\paragraph{Surrogates.}
\label{app:exp2-surrogates}

In addition to the same-PSF surrogate $\widetilde A$, we use the
Frobenius-optimal BCCB approximation
$\widetilde A_{\mathrm{opt}}=c(T)\otimes c(T)$, where $c(T)$ is
T.~Chan's optimal circulant approximation~\cite{chan1988optimal}.
Its first-column taps are $(1-|t|/N)g_t$, $N=256$, at periodic offsets
$t=-8,\ldots,8$. Thus $\widetilde A_{\mathrm{opt}}$ is circular
convolution with $p_{\mathrm{opt}}=g_{\mathrm{opt}}g_{\mathrm{opt}}^\top$
and has DC gain $0.98472$. For real $K$,
\[
\norm{T\otimes T-K\otimes K}_F^2
=\norm{T}_F^4+\norm{K}_F^4-2\langle T,K\rangle_F^2.
\]
This gives relative Frobenius distances of $0.10396$ for
$\widetilde A$ and $0.10322$ for $\widetilde A_{\mathrm{opt}}$.
The latter's signal-dependent mismatch is $0.07426$.

\paragraph{Cases.}
Five combinations of data, surrogate, and perturbation are evaluated:
\begin{center}
\small
\begin{tabular}{@{}llll@{}}
    \toprule
    Label & Data & Surrogate & Perturbation \\
    \midrule
    M-B & matched periodic $b_{\mathrm P}$ & $\widetilde A$ & bilinear \\
    Z-B & zero-boundary $b_0$ & $\widetilde A$ & bilinear \\
    M-P & matched periodic $b_{\mathrm P}$ & $\widetilde A$ & pure phase \\
    Z-P & zero-boundary $b_0$ & $\widetilde A$ & pure phase \\
    Z-C & zero-boundary $b_0$ & $\widetilde A_{\mathrm{opt}}$ & bilinear \\
    \bottomrule
\end{tabular}
\end{center}

Z-B is the primary transfer case and M-B its matched control. In each
row, the second model is the listed surrogate composed with $S$.
All spectral quantities, weights and thresholds use that row's
surrogate and data; the subscript in $D_b$ denotes $b_0$ or
$b_{\mathrm P}$. The targets are
$\tau\in\{0.05,0.075,0.10,0.125\}$, and all sensitivity
reconstructions use the roots $\lambda^\star$.

\section{Experiment 3: photographs, models, and fidelity}
\label{app:exp3-setup}

\paragraph{Preparing the photographs.}

The photographs are converted to grayscale using
$0.2989R+0.5870G+0.1140B$ on gamma-encoded values, without linearization.
The following square regions, specified as pixel boxes
(left, top, right, bottom) in the original files, are cropped:
\[
\begin{aligned}
\text{puffin:}\;&(0,0,2473,2473),\\
\text{owl:}\;&(600,0,3800,3200),\\
\text{goose:}\;&(0,300,2401,2701).
\end{aligned}
\]
Each crop is resized to $1110\times1110$ by Lanczos resampling and
stored as a 16-bit PNG with values in $[0,1]$. The full image is the
field $S$; its central $1024\times1024$ window is the ground truth
$x_{\mathrm{true}}$.

\paragraph{Data generation.}

The true PSF is the normalized disc of radius $43$, generated by MATLAB
\texttt{fspecial('disk',43)}. The noise-free data are
\[
y_c=\texttt{conv2}(S,\mathrm{psf},\texttt{'valid'}),
\]
computed by FFT. They have size $1024\times1024$ and are aligned with
the ground-truth window. For a standard Gaussian array $G$
(\texttt{twister}, seed $7$),
\[
y_0=y_c+\varepsilon\|y_c\|_F\,G/\|G\|_F,
\qquad \varepsilon\in\{0.03,0.01\}.
\]
The same $G$ is used at both noise levels. Data are neither clipped nor
quantized. Both model PSFs (radii $39$ and $43$) use the same
\texttt{fspecial} convention.

\paragraph{Original-model solves.}

PCG solves the normal equations for the truncated operator $A$ to
relative tolerance $10^{-8}$, with at most $500$ iterations and BCCB
preconditioner $(|\widehat{\mathrm{psf}}|^2+\lambda)^{-1}$.
Every solve meets the tolerance, requiring $44$--$223$ iterations
(\cref{tab:app-bird-fidelity}). No positivity constraint, clipping or
edge tapering is applied. Reconstructions live on an extended field:
$1102\times1102$ for radius $39$ and $1110\times1110$ for radius $43$.
All metrics use the native $1024\times1024$ window.

\paragraph{Model fidelity.}

Fidelity is measured on the inner window $W$, where each PSF footprint
lies within the known photograph, so no boundary assumption is needed.
The raw residual $\|b-Ax_{\mathrm{true}}\|_2/\|b\|_2$ is consistent
with the noise level for radius $43$. For radius $39$, subtracting noise
in quadrature gives an estimated model error of $1.2$--$1.7\%$ of the
data norm (\cref{tab:app-bird-fidelity}). A radius scan over $[30,60]$,
using an affine residual fit on a common inner window, selects the true
radius $43$ for every photograph and noise level.

\begin{table}[tb]
\centering
\caption{Experiment 3: model fidelity and operator mismatch. $r$: model disc radius. Raw and
affine residual: $\|b-Ax_{\mathrm{true}}\|_2/\|b\|_2$ on $W$ in percent, before and after the photometric fit
$b\approx c_1Ax_{\mathrm{true}}+c_0$. Radius: residual-minimizing disc radius of the scan. Mismatch: $\|(A-\widetilde A)x_{\mathrm{true}}\|_2/\|Ax_{\mathrm{true}}\|_2$
in percent (see below). PCG: range of iterations over all solves of that configuration.}
\label{tab:app-bird-fidelity}
\footnotesize\setlength{\tabcolsep}{3.6pt}
\begin{tabular}{@{}llrrrrrrrr@{}}
\toprule
 & noise & $r$ & raw & affine & $c_1$ & $c_0$ & radius & mismatch & PCG \\
\midrule
Puffin & 3\% & 39 & 3.29 & 3.26 & 0.987 & 0.008 & 43 & 0.87 & 49--140 \\
 & 1\% & 39 & 1.54 & 1.48 & 0.987 & 0.008 & 43 & 0.87 & 68--209 \\
 & 3\% & 43 & 3.10 & 3.10 & 1.000 & -0.000 & 43 & 1.21 & 50--147 \\
 & 1\% & 43 & 1.03 & 1.03 & 1.000 & -0.000 & 43 & 1.21 & 69--223 \\
\midrule
Owl & 3\% & 39 & 3.47 & 3.39 & 0.971 & 0.014 & 43 & 1.71 & 45--141 \\
 & 1\% & 39 & 1.98 & 1.84 & 0.971 & 0.014 & 43 & 1.71 & 55--209 \\
 & 3\% & 43 & 3.02 & 3.02 & 1.000 & -0.000 & 43 & 1.79 & 44--149 \\
 & 1\% & 43 & 1.01 & 1.01 & 1.000 & -0.000 & 43 & 1.79 & 51--222 \\
\midrule
Goose & 3\% & 39 & 3.24 & 3.21 & 0.989 & 0.004 & 43 & 1.88 & 45--139 \\
 & 1\% & 39 & 1.72 & 1.67 & 0.989 & 0.004 & 43 & 1.88 & 56--207 \\
 & 3\% & 43 & 2.91 & 2.91 & 1.000 & 0.000 & 43 & 1.84 & 46--145 \\
 & 1\% & 43 & 0.97 & 0.97 & 1.000 & 0.000 & 43 & 1.84 & 53--221 \\
\bottomrule
\end{tabular}
\end{table}

\paragraph{Errors against the original photograph.}

The main-text errors use the original photograph $x_{\mathrm{true}}$.
As a sensitivity check, we also use the affine reference
$c_1x_{\mathrm{true}}+c_0$ from \cref{tab:app-bird-fidelity}. Unit PSF
mass gives $A(c_1x+c_0)=c_1Ax+c_0$. For the radius-$43$ control,
$c_1=1$ and $c_0=0$ to three decimals. At radius $39$, the fit absorbs
some PSF mismatch ($c_1=0.971$--$0.989$) and reduces the oracle error
by $1.0\%$, $2.7\%$ and $1.0\%$ for puffin, owl and goose.

Selection uses $b$ alone and is unchanged by the reference. At
$\tau=0.015$ and $0.02$, changing the reference alters
$E/E_{\mathrm{or}}$ by at most $0.005$, without changing the main
conclusions. Some individual counts differ: at $\tau=0.05$, five
configurations are worse than the observation under the affine
reference, compared with four under the original reference.

The original-reference oracle is estimated by fitting a parabola through
the three stored errors around the minimum of the $61$ stored values.
It is not an evaluated reconstruction. This procedure reproduces the
stored scalar-search oracle error for the affine reference to relative
accuracy $1.5\times10^{-4}$. Oracle-associated targets are interpolated
from the stored root curve $\lambda^\star(\tau)$ in log--log
coordinates; the corresponding affine-reference targets are reproduced
to within $1\%$.

\paragraph{Operator mismatch.}

Because $A$ acts on an extended field and $\widetilde A$ on the native
window, $(A-\widetilde A)x_{\mathrm{true}}$ is not directly defined.
We compare them on a window $W$ obtained by removing a PSF-radius border
from the photograph. The original-model output there is
\[
\texttt{conv2}(x_{\mathrm{true}},\mathrm{psf},\texttt{'valid'}),
\]
while the comparison output is the reflexive blur of
$x_{\mathrm{true}}|_W$. They differ only within a PSF-radius band of
$\partial W$, where the reflexive extension replaces the actual scene
($15.8\%$ of the pixels of $W$ at radius $39$).

The relative mismatch is $0.9$--$1.9\%$
(\cref{tab:app-bird-fidelity}), compared with $7.5\%$ in Experiment~2.
Applying the zero-boundary/periodic comparison of Experiment~2 to these
photographs instead gives $0.123$, $0.113$ and $0.083$ for puffin, owl
and goose. Thus the smaller mismatch reflects the boundary treatment,
not simply the choice of photographs. A transfer model based on the
symmetric extension of the data would coincide with the surrogate on
the native window (mismatch about $5\times10^{-16}$); the truncated
model avoids that coincidence.

\section{Spectral diagnostics and conventions}
\label{app:exp-spectral-diagnostics}
For the half-pixel pairs of the experiments, $\mu_k=1$, and the thresholds
of \cref{thm:main-spectral-general,thm:phase-aware-threshold} are
\begin{equation}
    \lambda^\ddagger
    =\max_{k\in\mathcal K:\,a_k>s_k>0}\Bigl\{a_ks_k+(a_k+s_k)\sqrt{a_ks_k}\Bigr\},
    \qquad
    \lambda^\theta=\max_{k\in\mathcal K}\beta_k^\theta,
    \label{eq:exp-lambda-mon}
\end{equation}
with $\mathcal K=\{k:\widehat b_k\ne0\}$, $s_k=\rho_ka_k$, and
the maximum over an empty set equal to $0$. The thresholds depend on the surrogate, perturbation and active set,
but not on the data weights within that set. We treat all frequencies
as active; Gaussian noise makes every coefficient nonzero with
probability one.
We computed $\lambda^\theta$ only over the modes that violate
\eqref{eq:phase-dominance}. In Experiments~1 and~2 we took the largest
positive real root of each cubic $P_k$, which is $\beta_k^\theta$ of
\eqref{eq:beta-theta-def}. In Experiment~3 we took the largest positive
root at which $P_k$ changes sign, which is the supremum in
\eqref{eq:modal-sign-change}. The two values differ only at a terminal zero
where $P_k$ touches zero without changing sign. For the maximizing modes
reported below, the sign change holds, so both values agree there.
For the Gaussian surrogate, $\lambda^\ddagger=2.98859$ is attained at
the frequencies $(\pm1,0)$ and $(0,\pm1)$, which satisfy
\eqref{eq:phase-dominance} and therefore do not enter
$\lambda^\theta$. The value
$\lambda^\theta=2.97722$ is attained at $\pm(1,-1)$, where
$\theta_k=0$ (\cref{cor:half-pixel-thresholds}). For the Frobenius-optimal BCCB
surrogate of \cref{app:exp2-surrogates}, the corresponding values are
$2.89804$ and $2.88710$. With all frequencies active, these values are reproducible from the
kernels alone.

The bilinear multiplier and spatial action agree numerically to
$2.8\times10^{-16}$ and $4.0\times10^{-16}$, respectively.
For the pure-phase shift \eqref{eq:exp-purephase}, the discrepancy is the
directional roughness of \cref{prop:roughness}(i), including the removed
Nyquist components; it can be nonzero for noiseless data and is not an
error or noise estimate. On the $256\times256$ grid both perturbations have
$65{,}025$ nonzero multipliers, $99.61\%$ of which have nonzero phase.

Using the weights $w_k$ of \cref{app:fourier-implementation}, define
$\phi_k=a_k/(a_k^2+\lambda)$,
$\phi_k^S=s_k/(s_k^2+\lambda)$ and $h_k=\phi_k-\phi_k^S$.

\cref{prop:general-phase-formula,lem:phase-decomposition} lead to
\begin{equation}
    \Phi_b(\lambda)=\underbrace{\sum_k h_k^2\,w_k}_{\Phi_h(\lambda)}
    +\underbrace{\sum_k 2(1-\cos\theta_k)\,\phi_k\phi_k^S\,w_k}_{\Phi_{\mathrm{phase}}(\lambda)} .
    \label{eq:exp-decomposition}
\end{equation}
Primes denote differentiation with respect to $\lambda$ unless stated
otherwise. Define $[t]_+:=\max\{t,0\}$ and $[t]_-:=\min\{t,0\}$.
Differentiating term-wise, we can split the exact derivative as
$\Phi_b'=\mathcal R+\mathcal H+\mathcal P$, where
$\mathcal R=\sum_k[(h_k^2)']_+w_k\ge0$ collects modes in the increasing phase of
\cref{lem:scalar-filter-difference}, $\mathcal H=\sum_k[(h_k^2)']_-w_k\le0$ the decreasing modulus
terms, and $\mathcal P\le0$ the phase term, which is nonpositive mode by mode by
\cref{lem:phase-decomposition}. The \emph{rise ratio}
\begin{equation}
    r(\lambda)=\frac{\mathcal R(\lambda)}{|\mathcal H(\lambda)+\mathcal P(\lambda)|}
    \label{eq:exp-rise-ratio}
\end{equation}

satisfies $r(\lambda)<1$ exactly when $\Phi_b'(\lambda)<0$, provided
the denominator is nonzero. Otherwise, assess the sign directly from
$\mathcal R+\mathcal H+\mathcal P$. The rise ratio uses analytic
derivatives evaluated in floating point at sampled parameters; it gives
pointwise evidence, not a validated interval certificate. These spectral
identities apply to any nonzero $b$ satisfying the active-mode
hypotheses, whether or not $\widetilde A$ generated the data.

\paragraph{A check on a rectangular grid.}
\label{app:rect-grid}

On the square photographic grids, $\lambda^\theta=2.979$ is close to
$\lambda^\ddagger\approx3$. To examine the grid dependence identified
by \cref{thm:phase-aware-threshold}, we also use the Helsinki Deblur
Challenge 2021 geometry~\cite{juvonen_2021_4916176}: a $1460\times2360$ window,
mirror tile $2920\times4720$, and radius-$39$ disc PSF. This is an
operator calculation, without reconstruction. The bilinear pair gives
$\lambda^\theta=1.592$, $47\%$ below $\lambda^\ddagger=2.998$.

On this tile, $\omega_{1,k}+\omega_{2,k}=0$ requires
$k=j(73,-118)$, $j\in\mathbb Z$, since $\gcd(2920,4720)=40$. Already the
lowest such mode has $|\omega_k|=2\pi\sqrt2/40\approx0.22$~rad/px, beyond
the first zero of the disc spectrum ($\approx3.83/39\approx0.098$~rad/px),
where the disc strongly attenuates it. The bound of
\cref{cor:half-pixel-thresholds} is therefore far below $\lambda^\ddagger$.
The threshold is instead attained at the nearly anti-diagonal mode $k=(-13,21)$,
whose relative phase $\theta_k=9.1\times10^{-6}$ is nonzero but far too small for
\eqref{eq:phase-dominance}: the mode behaves almost as a phase-aligned one.
Its cubic changes sign at $\beta^\theta_k$, so by
\cref{thm:phase-aware-threshold}(iv) no data-independent threshold lies
below it. In either geometry, the experimentally selected weights lie well
below the thresholds; decrease in that range is verified using the
data-dependent certificate.

\section{Design of the paired noise-level study}
\label{app:tau-noise}
\label{sec:tau-noise}

This retrospective study supports \cref{tab:target-vs-parameter} and
\cref{sec:target-choice}. It examines how noise and boundary mismatch
affect the oracle-associated target, not how to select a target without
ground truth. Let $\lambda_{\mathrm{or}}^{\mathrm{ref}}$ be the best
parameter found by the search below, and set
$\tau_{\mathrm{or}}=D_b(\lambda_{\mathrm{or}}^{\mathrm{ref}})$.

\paragraph{Paired data.}
With the phantom and operators of \cref{app:exp2-solvers}, seeds
$7,\ldots,26$ generate one Gaussian array each. For a seed, we use the same
array at noise levels $0.003,0.01,0.03,0.05,0.10$, scaled by
each model's noise-free data norm. This gives $200$ positive-noise
datasets; one noiseless dataset per model brings the total to $202$.
The seed-$7$, $3\%$ case reproduces Experiment~2 to the reported
precision. All discrepancy evaluations use the same-PSF bilinear pair;
no interval certificates are computed.

\paragraph{Oracle search.}
We evaluate reconstruction error at $381$ logarithmically spaced weights
in $[10^{-16},10^3]$ and refine every sampled interior local minimum
using \texttt{fminbnd} in $\log\lambda$ (\texttt{TolX}$=10^{-9}$).
The best refined or grid value, including endpoints, defines
$\lambda_{\mathrm{or}}^{\mathrm{ref}}$ and
$E_{\mathrm{or}}^{\mathrm{ref}}$. This numerical oracle is not a
certified global minimum. It uses a broader search than the $61$-point
$[10^{-6},1]$ oracle in the main experiments and does not change their
reported results. Original-model errors use the eigendecomposition of
\cref{app:exp2-solvers}, without clipping, positivity constraints or
spectral cutoffs. All $200$ positive-noise optima are interior. The two
noiseless optima occur at the lower bound and are excluded from the
summary.

\begin{table}[!ht]
\centering
\caption{Paired noise-level study: medians over $20$ realizations. The
oracle-associated targets use ground truth. The two noiseless cases are
excluded because their searched optima occur at the lower search boundary.}
\label{tab:tau-noise-medians}
\small
\begin{tabular}{@{}rrrrr@{}}
\toprule
& \multicolumn{2}{c}{Matched periodic} & \multicolumn{2}{c}{Zero-boundary transfer}\\
\cmidrule(lr){2-3}\cmidrule(lr){4-5}
Noise (\%) & $\lambda_{\mathrm{or}}^{\mathrm{ref}}$ & $\tau_{\mathrm{or}}$
& $\lambda_{\mathrm{or}}^{\mathrm{ref}}$ & $\tau_{\mathrm{or}}$\\
\midrule
0.3 & 0.000115 & 0.07983 & 0.000139 & 0.43709 \\
1 & 0.000883 & 0.07357 & 0.00102 & 0.19426 \\
3 & 0.00514 & 0.06717 & 0.00558 & 0.11064 \\
5 & 0.0107 & 0.06429 & 0.0112 & 0.09206 \\
10 & 0.0253 & 0.06120 & 0.0254 & 0.07719 \\
\bottomrule
\end{tabular}
\end{table}

\paragraph{Oracle parameters and targets.}
The oracle-associated target decreases with noise in every paired
realization of both models (\cref{tab:tau-noise-medians}), but its level
depends on the boundary model. At $0.3\%$ noise, the median targets are
$0.07983$ and $0.43709$; at $10\%$, they are $0.06120$ and $0.07719$.
The corresponding oracle parameters are much closer across the two
models. Thus similar desirable regularization levels need not imply
similar surrogate-discrepancy targets.

% ------------------------------------------------------------
% Supplement (merged from supplement.tex; sits after the bibliography)
% ------------------------------------------------------------
% \clearpage
% \input{supplement}

\end{document}